\documentclass[english,12pt]{amsart}
\usepackage{amssymb}
\usepackage{amsfonts}
\usepackage[foot]{amsaddr}
\usepackage{amscd}
\usepackage[all]{xy}
\usepackage{color}
\usepackage{enumerate}
\usepackage[utf8]{inputenc}

\usepackage[T1]{fontenc}
\usepackage{xr-hyper}
\usepackage[colorlinks=true, breaklinks=true, urlcolor= black, allcolors= black,bookmarksopen=true,linktocpage=true,plainpages=false,pdfpagelabels]{hyperref}
\usepackage{cleveref}

\CompileMatrices

\swapnumbers

\def\rad{\operatorname{rad}}
\def\radop{\rad_{\mathrm{op}}}
\def\radcl{\rad_{\mathrm{cl}}}

\def\radopeqc{\rad_{\mathrm{op},\eqc}}

\def\radcleqc{\rad_{\mathrm{cl},\eqc}}
\def\mix{{\mathrm{mix}}}
\def\heqc{h$^\eqc$}
\def\hc{h$^{{\rm{coars}}}$}
\def\hmix{h$^\mix$}

\def\val{{\mathrm{val}}}
\def\Lval{\cL_\val}

\def\eqc{{\mathrm{ecc}}}  

\def\11{{\mathbf 1}}

\def\CC{{\mathbb C}}

\def\FF{{\mathbb F}}

\def\NN{{\mathbb N}}

\def\QQ{{\mathbb Q}}
\def\RR{{\mathbb R}}

\def\ZZ{{\mathbb Z}}

\def\cA{{\mathcal A}}

\def\cL{{\mathcal L}}
\def\cM{{\mathcal M}}

\def\cO{{\mathcal O}}

\def\cT{{\mathcal T}}

\def\llp{\mathopen{(\!(}}

\def\rrp{\mathopen{)\!)}}

\mathchardef\alphag="7C0B \mathchardef\betag="7C0C
\mathchardef\gammag="7C0D \mathchardef\deltag="7C0E
\mathchardef\varepsilong="7C22 \mathchardef\varphig="7C27
\mathchardef\psig="7C20 \mathchardef\zetag="7C10
\mathchardef\epsilong="7C0F \mathchardef\rhog="7C1A
\mathchardef\taug="7C1C \mathchardef\upsilong="7C1D
\mathchardef\iotag="7C13 \mathchardef\thetag="7C12
\mathchardef\pig="7C19 \mathchardef\sigmag="7C1B
\mathchardef\etag="7C11 \mathchardef\omegag="7C21
\mathchardef\kappag="7C14 \mathchardef\lambdag="7C15
\mathchardef\mug="7C16 \mathchardef\xig="7C18
\mathchardef\chig="7C1F \mathchardef\nug="7C17
\mathchardef\varthetag="7C23 \mathchardef\varpig="7C24
\mathchardef\varrhog="7C25 \mathchardef\varsigmag="7C26
\mathchardef\Omegag="7C0A \mathchardef\Thetag="7C02
\mathchardef\Sigmag="7C06 \mathchardef\Deltag="7C01
\mathchardef\Phig="7C08 \mathchardef\Gammag="7C00
\mathchardef\Psig="7C09 \mathchardef\Lambdag="7C03
\mathchardef\Xig="7C04 \mathchardef\Pig="7C05
\mathchardef\Upsilong="7C07

\newtheorem{thm}[subsubsection]{Theorem}
\newtheorem{lem}[subsubsection]{Lemma}
\newtheorem{cor}[subsubsection]{Corollary}
\newtheorem{prop}[subsubsection]{Proposition}

\newtheorem{addendum}
{Addendum}
\Crefname{prop}{Proposition}{Propositions}
\Crefname{thm}{Theorem}{Theorems}
\Crefname{lem}{Lemma}{Lemmas}
\Crefname{cor}{Corollary}{Corollaries}
\Crefname{conj}{Conjecture}{Conjectures}

\theoremstyle{definition}
\newtheorem{defn}[subsubsection]{Definition}

\newtheorem{notn}[subsubsection]{Notation}
\newtheorem{example}[subsubsection]{Example}

\newtheorem{def-prop}[subsubsection]{Proposition-Definition}
\newtheorem{def-theorem}[subsubsection]{Theorem-Definition}
\newtheorem{def-lem}[subsubsection]{Lemma-Definition}

\Crefname{defn}{Definition}{Definitions}

\theoremstyle{remark}
\newtheorem{remark}[subsubsection]{Remark}

\newtheorem{question}[subsubsection]{Question}

\theoremstyle{plain}

\numberwithin{equation}{subsection}

\def\boxit#1#2{\setbox1=\hbox{\kern#1{#2}\kern#1}%
\dimen1=\ht1 \advance\dimen1 by #1 \dimen2=\dp1 \advance\dimen2 by
#1
\setbox1=\hbox{\vrule height\dimen1 depth\dimen2\box1\vrule}%
\setbox1=\vbox{\hrule\box1\hrule}%
\advance\dimen1 by .4pt \ht1=\dimen1 \advance\dimen2 by .4pt
\dp1=\dimen2 \box1\relax}

\renewcommand{\theequation}{\thesubsection.\arabic{equation}}

\mathchardef\alphag="7C0B \mathchardef\betag="7C0C
\mathchardef\gammag="7C0D \mathchardef\deltag="7C0E
\mathchardef\varepsilong="7C22 \mathchardef\varphig="7C27
\mathchardef\psig="7C20 \mathchardef\zetag="7C10
\mathchardef\epsilong="7C0F \mathchardef\rhog="7C1A
\mathchardef\taug="7C1C \mathchardef\upsilong="7C1D
\mathchardef\iotag="7C13 \mathchardef\thetag="7C12
\mathchardef\pig="7C19 \mathchardef\sigmag="7C1B
\mathchardef\etag="7C11 \mathchardef\omegag="7C21
\mathchardef\kappag="7C14 \mathchardef\lambdag="7C15
\mathchardef\mug="7C16 \mathchardef\xig="7C18
\mathchardef\chig="7C1F \mathchardef\nug="7C17
\mathchardef\varthetag="7C23 \mathchardef\varpig="7C24
\mathchardef\varrhog="7C25 \mathchardef\varsigmag="7C26
\mathchardef\Omegag="7C0A \mathchardef\Thetag="7C02
\mathchardef\Sigmag="7C06 \mathchardef\Deltag="7C01
\mathchardef\Phig="7C08 \mathchardef\Gammag="7C00
\mathchardef\Psig="7C09 \mathchardef\Lambdag="7C03
\mathchardef\Xig="7C04 \mathchardef\Pig="7C05
\mathchardef\Upsilong="7C07

\newcommand{\RV}{\mathrm{RV}}

\newcommand{\rv}{\operatorname{rv}}

\newcommand{\Th}{\operatorname{Th}}
\newcommand{\valring}{\mathcal{O}_K}

\newcommand\LI{\cL_I}

\newcommand{\eq}{^{\mathrm{eq}}}

\newcommand{\grad}{\operatorname{grad}}

\newcommand{\abs}[1]{\lvert#1\rvert}

\definecolor{immi}{rgb}{0,.6,.1}

\newbox\removebox
\newcommand\remove[1]{%
\setbox\removebox=\ifmmode\hbox{$#1$}\else\hbox{#1}\fi%
\leavevmode
\rlap{\textcolor{blue}{\vrule height0.8ex depth-0.6ex width\wd\removebox}}%
\box\removebox
}
\long\def\bigremove#1{%
\par\setbox\removebox=\vbox{#1}%
\vbox{%
\vbox to0pt{\hbox{\tikz\draw[color=blue,thick] (0,0) -- (\wd\removebox,-\ht\removebox)  (\wd\removebox,0) -- (0,-\ht\removebox);}}
\box\removebox
}
}

\newcommand\tp{\mathrm{tp}}
\newcommand\qftp{\mathrm{qftp}}
\newcommand\acl{\mathrm{acl}}
\newcommand\dcl{\mathrm{dcl}}
\newcommand{\restr}[2]{{\left.#1\right|_{#2}}}

\definecolor{orange}{rgb}{1,0.5,0}

\newcommand{\private}[1]{\leavevmode{\scriptsize\color{blue}\marginpar{{\scriptsize Private comment}}#1\par}}
\renewcommand{\private}[1]{}

\thanks{The author R.\,C. was partially supported by the European Research Council under the European Community's Seventh Framework Programme (FP7/2007-2013) with ERC Grant Agreement nr.\ 615722 MOTMELSUM, by KU Leuven IF C14/17/083, and thanks the Labex CEMPI  (ANR-11-LABX-0007-01). The author I.\,H. was partially supported by the research training group
\emph{GRK 2240: Algebro-Geometric Methods in Algebra, Arithmetic and Topology} and by the individual research grant No.~426488848, both funded by the Deutsche Forschungsgemeinschaft (DFG, German Research Foundation).
The author S.\,R. was partially supported by GeoMod AAPG2019 (ANR-DFG), Geometric and Combinatorial Configurations in Model Theory.
The author F.\,V. was partially supported by KU Leuven IF C14/17/083, and partially by F.W.O. Flanders (Belgium) with grant number 11F1921N}

\title[Hensel minimality II]{Hensel minimality II: Mixed characteristic and a diophantine application}

\author[Cluckers]{Raf Cluckers}
\address{Raf Cluckers, Univ.~Lille, CNRS, UMR 8524 - Laboratoire Paul Painlev\'e, F-59000 Lille, France, and,
KU Leuven, Department of Mathematics, B-3001 Leu\-ven, Bel\-gium}
\email{Raf.Cluckers@univ-lille.fr}
\urladdr{http://rcluckers.perso.math.cnrs.fr/}

\author[Halupczok]{Immanuel Halupczok}
\address{Immanuel Halupczok, Lehrstuhl f\"ur Algebra und Zahlentheorie, Mathematisches Institut, Universit\"atsstr. 1, 40225 D\"usseldorf, Germany}
\email{math@karimmi.de}
\urladdr{http://www.immi.karimmi.de/en/}

\author[Rideau]{Silvain Rideau-Kikuchi}
\address{Silvain Rideau-Kikuchi, CNRS, Département des mathématiques et applications, École normale supérieure, 45 rue d’Ulm, F-75230 Paris Cedex 05}
\email{silvain.rideau@imj-prg.fr}
\urladdr{https://webusers.imj-prg.fr/\~{}silvain.rideau/en/}

\author[Vermeulen]{Floris Vermeulen}
\address{Floris Vermeulen, KU Leuven, Department of Mathematics, B-3001 Leu\-ven, Bel\-gium}
\email{floris.vermeulen@kuleuven.be}
\urladdr{https://sites.google.com/view/floris-vermeulen/homepage}



\begin{document}

\begin{abstract}
In this paper together with the preceding Part I \cite{CHR}, we develop a framework for tame geometry on Henselian valued fields of characteristic zero, called Hensel minimality.
It adds to \cite{CHR} the treatment of the mixed characteristic case.
Hensel minimality is inspired by o-minimality and its role in real geometry and diophantine applications.
We develop geometric results and applications for Hensel minimal structures that were previously known only under stronger or less axiomatic assumptions, and which often have counterparts in o-minimal structures.
We prove a Jacobian property, a strong form of Taylor approximations of definable functions, resplendency results and cell decomposition, all under Hensel minimality, more precisely, $1$-h-minimality.
We obtain a diophantine application of counting rational points of bounded height on Hensel minimal curves.
\end{abstract}

\maketitle

\newpage

\setcounter{tocdepth}{1}
\tableofcontents

\section{Introduction}
\label{sec:intro}

\subsection{}\label{ssec:intr}

For a long time one has been looking for adequate analogues of o-minimality in other settings than for real geometry. Several notions  have been put forward, each with certain strengths and weaknesses. In Part I \cite{CHR} together with this sequel paper, we put forward a notion for tame geometry on non-archimedean valued fields, called Hensel minimality, both simple and strong and providing a common variant for the settings from, among other, \cite{Haskell,Macpherson,HM,HK,CLb}.

\par

Let us right away explain the key flavor of the definition of Hensel minimality.
As often, finding a correct analogue involves some reformulations. 
Here is a way to reformulate o-minimality: a structure on the real field $\RR$ is o-minimal if and only if for each definable subset $X$ of $\RR$, there exists a finite set $C\subset \RR$, such that, for any $x$ in $\RR$, the condition $x\in X$ only depends on the signs of $x-c$ for $c$ in $C$, the sign being negative, zero, or positive. In an equicharacteristic zero valued field $K$, we replace the sign map by the projection map $\rv:K\to K/(1+\cM_K)$, where $\cM_K$ is the maximal ideal of the valuation ring of $K$, and $K/(1+\cM_K)$ the quotient of multiplicative semi-groups. In the mixed characteristic case, one additionally uses ideals of the form $N\cdot\cM_K$ for nonzero integers $N$, instead of just $\cM_K$ itself. Apart from replacing the sign map, we keep the condition of every definable subset of $K$ being ``controlled'' by a finite set $C\subset K$. In the real case, $C$ is automatically definable over the same parameters as $X$; in the valued field case, this is no longer automatic and needs to be imposed. Depending on the precise kind of parameters one allows, one obtains
the notion of $0$-h-minimality, $\omega$-h-minimality, or something in-between: $\ell$-h-minimality for integers $\ell>0$; see \Cref{sec:mixed} for the detailed definitions. When we do not want to be precise about the specific version, we just say ``Hensel minimality''.

\par

Hensel minimality is similar to o-minimality \cite{VdD,PillaySteinhornI} not only in its definition, but also for having strong consequences. In Part~I \cite{CHR} we focused on the equicharacteristic zero case. In this sequel, we focus on the mixed characteristic case. We obtain analogues of many of the results from \cite{CHR}. In addition, in \Cref{sec:counting}, we give a new diophantine application, similar to the results by Pila-Wilkie \cite{PW} from in o-minimal case; for the moment this application is for curves only. Many proofs from \cite{CHR} adapt to mixed characteristic in a rather straightforward way; nevertheless, many proofs are repeated in this sequel, and for many others, we sketch the key ideas, so that the reader has the choice between reading the short version here or the detailed version in \cite{CHR}.

There are several ways to adapt the notions of Hensel minimality from \cite{CHR} to
the mixed characteristic case, based on coarsenings, compactness, or on a more literal adapation of the definitions for the equi-characteristic zero case. We treat several of these variants, and in one of the most important cases (namely $1$-h-minimality) we show that they are equivalent. From this equivalence we derive that many of the geometric results of \cite{CHR} also hold in the mixed characteristic case, with in particular some strong results on Taylor approximation of definable functions. In \cite{Verm:h-min}, equivalences of these notions are shown also in other cases.

The name ``Hensel minimality'' comes from the intuition that a Hensel minimal valued field behaves as nicely as a henselian valued field in the pure valued field language.
In equi-characteristic $0$, this is reflected by the fact that a pure valued field $K$ is Hensel minimal if and only if $K$ is henselian (where any version of Hensel minimality can be used). In mixed characteristic, this is still true if $K$ is finitely ramified (see Remark~\ref{rem:hensel}), but in general, Hensel minimality only captures tameness from the ramification degree on.



\subsection{}\label{ssec:intr.b}
Let us give a short overview of the paper.
In \Cref{sec:mixed} we introduce several variants of Hensel minimality (for both the mixed and the equi-characteristic cases), and we prove that all of them are equivalent in the key case of $1$-h-minimality (\Cref{thm:equivalences:1hmin}). We do this for a large part by classical strategies: coarsening of valuations of mixed characteristic to equi-characteristic zero, and model theoretic compactness. A combination of these classical strategies with new geometric and model theoretic arguments lead to our strongest results, like \Cref{cor:high-ord:mix} on Taylor approximation. It is precisely this result \Cref{cor:high-ord:mix} on Taylor approximation that plays a key role in our diophantine application \Cref{thm:QQp-c-dim}, where we estimate the number of rational points of bounded height on transcendental curves. This is the Hensel minimal analogue of point counting on transcendental curves case as treated by Bombieri and Pila \cite{Bombieri-Pila} and o-minimally by Pila-Wilkie \cite{PW}, and, it is the axiomatic analogue of point counting on subanalytic sets as studied in \cite{CCL-PW} and \cite{CFL}, and on analytic sets in \cite{BinKat}.

In \Cref{sec:respl}, we develop some resplendency results analogous to the ones from 
\cite[Section~4]{CHR}, namely, that Hensel minimality is preserved under certain expansions of the structure (e.g.\ arbitrary expansions of the value group and the residue field).
Those results are used to show the equivalences of \Cref{thm:equivalences:1hmin}.
It follows from the equivalence with item (4) in \Cref{thm:equivalences:1hmin} that $1$-h-minimality is preserved under coarsening of the valuation; such a coarsening result was previously known only for $\omega$-h-minimality by 
\cite[Corollary~4.2.4]{CHR}. Note that in \cite{Verm:h-min}, several results of this paper and of part I \cite{CHR} are extended from $1$-h-minimality to $\ell$-h-minimality for each $\ell\ge 1$, for example, their preservation under coarsenings of the valuation.

We end the paper with some open questions.
One of the big challenges in the current framework is to push our diophantine application further towards arbitrary dimension, and thus to get a full, Hensel minimal analogue of the results by Pila-Wilkie \cite{PW}.
Other main challenges are further developments of the geometry, like t-stratifications building on and extending 
\cite[Section~5.5]{CHR}.  It may also be interesting to compare the notions of Hensel minimal structures with distal expansions on valued fields, as studied recently in \cite{ACGZ}.

\section{Hensel minimality in characteristic zero (including mixed characteristic)}\label{sec:mixed}

In this whole section we fix a theory $\cT$ of valued fields of characteristic zero in a language $\cL$ expanding the language $\Lval= \{+,\cdot,\cO_K\}$ of valued fields. Note that $\cT$ is allowed to be non-complete, and that each model $K$ of $\cT$ is a valued field of characteristic zero, which includes both possibilities that  $K$ has mixed characteristic or equicharacteristic zero.

In this section we give four alternative definitions of $1$-h-minimality for $\cT$ (see Definitions \ref{defn:mixed:RVlambdaRVm}, \ref{defn:mixed:eqc}, \ref{defn:mixed:c},  \ref{defn:T1T2}), and we show that they are equivalent in \Cref{thm:equivalences:1hmin}. To keep the generality of \cite{CHR}, we will treat more generally $\ell$-h-minimality for $\ell\ge 0 $ either an integer or equal to $\omega$. The first definition is a close adaptation of the main notion of Hensel minimality of 
\cite[Definition~2.3.3]{CHR}.  Definitions \ref{defn:mixed:eqc} and \ref{defn:mixed:c} are based on coarsenings, and \Cref{defn:T1T2} corresponds to the criterion for $1$-h-minimality from 
\cite[Theorem~2.9.1]{CHR}. Since for $\ell=1$ these notions coincide, we simply will call them $1$-h-minimality.

\subsection{Basic terminology}\label{ssec:notation}

We use the following terminology, notation and concepts from \cite{CHR}.
By a valued field we mean a non-trivially valued field, i.e., the field of fractions of a valuation ring which is not a field itself. Any valued field $K$ is a structure in the language $\Lval= \{+,\cdot,\cO_K\}$ of valued fields, where $+$ and $\cdot$ are addition and multiplication on $K$, and where $\cO_K$ is (a predicate for) the valuation ring of $K$. The maximal ideal of $\cO_K$ is denoted by $\cM_K$.
We use multiplicative notation for the value group which we denote by $\Gamma_K^\times$, and we write  $\Gamma_K :=  \Gamma_K^\times \cup\{0\}$. We write $|\cdot|\colon K \to \Gamma_K$ for the valuation map.
By an open ball we mean a set of the form $B_{<\lambda}(a)=\{x\in K\mid |x-a| < \lambda\}$ for some $\lambda$ in $\Gamma_K^\times$ and some $a$ in $K$. We define $\radop B_{<\lambda}(a)=\lambda$ to be the open radius of such a ball. Similarly, a closed ball is a set of the form $B_{\leq \lambda}(a) = \{x\in K\mid |x-a|\leq \lambda\}$. Its closed radius is $\radcl B_{\leq \lambda}(a)=\lambda$.

For any proper ideal $I$ of the valuation ring $\cO_K$, we write $\RV_I$ for the corresponding leading term structure, i.e., the disjoint union of $\{0\}$ with the quotient group $K^\times / (1+I)$, and $\rv_I:K\to \RV_I$ for the leading term map, i.e., the quotient map extended by sending sending $0$ to $0$.
When $I$ is the open ball $\{x \in K \mid |x| < \lambda\}$ for some  $\lambda \leq 1$ in $\Gamma_K^\times$, we simply write $\RV_\lambda$ and $\rv_\lambda$ instead of $\RV_I$ and $\rv_I$, and we write $\RV$ and $\rv$ instead of $\RV_1$ and $\rv_1$.

We will sometimes also write $\RV_K$ for $\RV$ and similarly $\RV_{K,\lambda}$ for $\RV_\lambda$ if multiple fields are under consideration.

\begin{defn}
Let $\lambda \leq 1$ be an element of $\Gamma_K^\times$.
\begin{enumerate}
 \item Given an arbitrary set $X \subset K$ and a finite non-empty set $C \subset K$,
 we say that $X$ is
 $\lambda$-prepared by $C$ if the condition whether some $x\in K$ lies in $X$ depends only on the tuple $(\rv_\lambda(x-c))_{c\in C}$.
\item
We say that a ball $B \subset K$ is $\lambda$-next to an element $c \in K$ if
$$
B= \{x\in K\mid \rv_\lambda(x-c) = \xi \}
$$
for some (nonzero) element $\xi$ of $\RV_\lambda$.
\item
We say that a ball $B \subset K$ is $\lambda$-next to a finite non-empty set $C \subset K$ if $B$ equals $\bigcap_{c\in C} B_c$ with $B_c$ a ball $\lambda$-next to $c$ for each $c\in C$.
\end{enumerate}
\end{defn}

One easily verifies that for fixed $C$ and $\lambda$, the set of all balls $\lambda$-next to $C$ forms a partition of $K \setminus C$. Moreover, a set $X \subset K$ is $\lambda$-prepared by $C$ if and only if $X \setminus C$ is a union of parts of this partition (without any condition on $X \cap C$).
Decreasing $\lambda$ refines this partition, i.e., if $\lambda' < \lambda$ and $C$ $\lambda$-prepares $X$, then it also $\lambda'$-prepares $X$. Similarly, if $C' \supset C$ and $C$ $\lambda$-prepares $X$, then $C'$ also $\lambda$-prepares $X$.

\begin{example}
Fix any $\lambda \le 1$ in $\Gamma_K^\times$.
\begin{enumerate}
 \item
A finite set $X \subset K$ is $\lambda$-prepared by $C$ if and only if $C$ contains $X$.
 \item
A subset $X \subset K$ is $\lambda$-prepared by the set $C = \{0\}$ if and only if $X$ is a (possibly infinite) union of fibers of the map $\rv_\lambda$, i.e., if it is of the form $X = \rv_\lambda^{-1}(\Xi)$, for an arbitrary subset $\Xi \subset \RV_\lambda$.
\item Every open ball $B = B_{< \lambda}(a)$ of radius $\lambda$ contained in $\valring$ is $\lambda$-prepared by $C = \{0\}$, though $B$ might not be $\lambda$-next to $C$. (It is $\lambda$-next to $C$ if and only if $|a| = 1$.)
\end{enumerate}
\end{example}

The following notation specific to mixed characteristic already appears in 
\cite[Section~6]{CHR}.
\begin{defn}[Equi-characteristic $0$ coarsening]\label{not.eqc}
Given a model $K \models \cT$, we write $\cO_{K,\eqc}$ for the smallest subring of $K$ containing $\cO_K$ and $\QQ$ and we let $|\cdot|_{\eqc}\colon K \to \Gamma_{K,\eqc}$ be the corresponding valuation. (Thus, $|\cdot|_{\eqc}$ is the finest coarsening of $|\cdot|$ which has equi-characteristic $0$; note that $|\cdot|_{\eqc}$ can be a trivial valuation on $K$.) If $|\cdot|_{\eqc}$ is a nontrivial valuation (i.e., $\cO_{K,\eqc} \ne K$), then we also use the following notation: $\rv_{\eqc}\colon K \to \RV_\eqc$ is the leading term map with respect to $|\cdot|_{\eqc}$; given $\lambda \in \Gamma_{K,\eqc}$, $\rv_{\lambda}\colon K \to \RV_\lambda$ is the leading term map with respect to $\lambda$; and $\cL_{\eqc}$ is the extension of $\cL$ by a predicate for $\cO_{K,\eqc}$. More generally, for any non-trivial coarsening $|\cdot|_c$ of the valuation on $K$, write $\cL_c$ for the extension of $\cL$ by a predicate for the valuation ring for $|\cdot|_c$.
\end{defn}

\subsection{Equivalent definitions}\label{ssec:mixed:RVlambda}

There are several natural notions of  Hensel minimality in mixed characteristic. We give four possible definitions. \Cref{thm:equivalences:1hmin} states that these are all equivalent, in the case of $1$-h-minimality.

We first adapt the main definition of Hensel minimality from 
\cite[Definition~2.3.3]{CHR} to include the mixed characteristic case. Recall that $\cT$ is a theory of valued fields of characteristic zero in  a language $\cL\supset \Lval$.

\begin{defn}\label{defn:mixed:RVlambdaRVm}
Let $\ell\geq 0$ be either an integer or $\omega$. Say that $\cT$ is \emph{$\ell$-\hmix-minimal} if for each model $K$ of $\cT$, for each integer $n \ge 1$, each $\lambda \le 1$ in $\Gamma_K^\times$, each $A\subset K$, each finite $A'\subset \RV_\lambda$ of cardinality $\# A'\le \ell$, and each $(A\cup \RV_{|n|} \cup A')$-definable set $X\subset K$, there exists an integer $m \ge 1$ such that $X$ is $|m|\lambda$-prepared by a finite $A$-definable set $C\subset K$.
If all models of $\cT$ are of equicharacteristic zero,
we also call $\cT$ simply \emph{$\ell$-h-minimal}.
\end{defn}

As one sees from the definition,
bigger $\ell$ yield stronger conditions, i.e., for $\ell < \ell'$, $\ell'$-\hmix-minimality implies $\ell$-\hmix-minimality.

Clearly, if all models $K$ of $\cT$ are of equicharacteristic zero, then one can take $n=m=1$, so
the above definition of $\ell$-h-minimality agrees with 
\cite[Definition~2.3.3]{CHR}.
In the mixed characteristic case, one is obliged to take the valuation of integers into account in \Cref{defn:mixed:RVlambdaRVm}. The general philosophy is that $\lambda$-preparation in equi-characteristic $0$ becomes $\lambda\cdot |m|$-preparation  in mixed characteristic, for some integer $m \ge 1$, as the following example illustrates.

\begin{example}\label{eg:cubes.Q3}
The set $X$ of cubes in the $3$-adic numbers $\QQ_3$ cannot be $1$-prepared by any finite set $C$, since each of the infinitely many disjoint balls $27^r(1+3\ZZ_3)$, $r \in \ZZ$, contains both, cubes and non-cubes. However, $X$ is a union of (infinitely many) fibers of the map $\rv_{|3|}\colon \QQ_3 \to \RV_{|3|}$, so it is $|3|$-prepared by the set $\{0\}$.
\end{example}

In \cite{DG.tame}, Dolich and Goodrick introduced the notion of a \emph{visceral} structure (a structure equipped with a uniformly definable topological base satisfying certain simple axioms) and showed several tameness results in this abstract setting. It is easy to see that $0$-\hmix-minimality implies viscerality in this sense.

\begin{remark}\label{rem:hensel}
In the equi-characteristic $0$ case, $0$-h-minimality implies definable spherical completeness (meaning that every definable family of nested balls has non-empty intersection), and this in turn implies henselianity (see \cite{CHR}, Lemma~2.7.1 and Theorem~2.7.2). A similar implication holds in mixed characteristic if we restrict to finitely ramified fields $K$: If the theory such a $K$ is $0$-\hmix-minimal, then it is definable spherical complete by \cite[Proposition~1.4]{iBW.sph}, and this in turn implies that $K$ is henselian using the usual argument via Newton approximation (see the proof of \cite[Theorem~2.7.2]{CHR} for details).
In contrast, note that if $K$ is not finitely ramified, then $0$-\hmix-minimality does not imply definable spherical completeness; see Example~1.5 in \cite{iBW.sph}.
\end{remark}

The following two definitions use equicharacteristic zero coarsenings of the valuation to define mixed characteristic Hensel minimality in terms of the equicharacteristic notions from 
\cite[Definition~2.3.3]{CHR} (which we just recalled in Definition~\ref{defn:mixed:RVlambdaRVm}).
We use notation from \Cref{not.eqc}.

\begin{defn}
\label{defn:mixed:eqc}
Let $\ell\geq 0$ be either an integer or $\omega$. We say that $\cT$ is \emph{$\ell$-\heqc-minimal} if for every model $K \models \cT$ the following holds: If the valuation $|\cdot|_{\eqc}$ on $K$ is non-trivial, then the $\cL_{\eqc}$-theory of $K$, when considered as a valued field with the valuation $|\cdot|_{\eqc}$, is $\ell$-h-minimal.
\end{defn}

One can also require every equi-characteristic zero coarsening to be $\ell$-h-minimal, leading to the following definition.

\begin{defn}
\label{defn:mixed:c}
Let $\ell\geq 0$ be either an integer or $\omega$. We say that $\cT$ is \emph{$\ell$-\hc-minimal} if for every model $K \models \cT$ and each non-trivial equicharacteristic coarsening $|\cdot|_c$ of the valuation on $K$, the $\cL_{c}$-theory of $K$, when considered as a valued field with the valuation $|\cdot|_{c}$, is $\ell$-h-minimal.
\end{defn}

\begin{remark}\label{rem:equichar}
Since in an equi-characteristic zero valued field $K$, we have $\cO_{K,\eqc} = \cO_K$,
the theory of such a field $\ell$-\heqc-minimal if and only if it is $\ell$-h-minimal, for each $\ell\ge 0$.
For more subtle reasons, $\omega$-\hc-minimality is equivalent to $\omega$-\heqc-minimality by 
\cite[Corollary~4.2.4]{CHR}. Even more, in \cite{Verm:h-min} it is shown that in equi-characteristic zero, $\ell$-\hc-minimality is equivalent to $\ell$-\heqc-minimality for any $\ell\geq 1$. In mixed characteristic, we will show the analogue of this (and more) for $\ell=1$ in Theorem \ref{thm:equivalences:1hmin}.
\end{remark}

By the usual play of compactness, preparation results that hold in each model of $\cT$ also hold uniformly for all models of $\cT$, and results in equi-characteristic zero can be transfered to mixed characteristic.
We will give examples in the proofs for \Cref{cor:prep:nc} and \Cref{cor:prep:mix} to illustrate how compactness is used for these purposes.

In 
\cite[Theorem~2.9.1]{CHR}, a geometric criterion for $1$-h-minimality is given in equi-characteristic zero. The following is a mixed characteristic version of that criterion.

\begin{defn}
\label{defn:T1T2}
Let $f: K\to K$ be $A$-definable, for some set $A\subset
 K\cup \RV_{|n|}$ for some positive integer $n$. We define the following two properties:
\begin{itemize}
\item[(T1$^\mix$)] There exists a finite $A$-definable $C\subset   K$ and a positive integer $m$ such that for every ball $B$ $|m|$-next to $C$ there exists $\mu_B\in \Gamma_K$ such that for $x,y\in B$ we have
\[
\mu_B\cdot|m|\cdot |x-y|\leq  |f(x)-f(y)| \leq \mu_B\cdot \left|\frac{1}{m}\right| \cdot  |x-y|.
\]
\item[(T2)] The set $\{y\in K\mid f^{-1}(y) \text{ is infinite}\}$ is finite.
\end{itemize}
We say that $\cT$ satisfies (T1,T2) if for all $f$, $A$ and $n$ as above, the two conditions (T1$^\mix$) and (T2) hold.
\end{defn}

We now have several variants of 1-h-minimality in mixed characteristic. The main result of this section is the following, stating that all of the above definitions agree.

\begin{thm}\label{thm:equivalences:1hmin}
The following are equivalent, for a theory $\cT$ of valued fields of characteristic zero (possibly of mixed characteristic) in a language $\cL$ containing $\Lval$.
\begin{enumerate}
\item $\cT$ is 1-\heqc-minimal.
\item $\cT$ is 1-\hmix-minimal.
\item $\cT$ satisfies (T1,T2).
\item $\cT$ is 1-\hc-minimal.
\end{enumerate}
Therefore we will call this common notion simply \emph{1-h-minimality}.
\end{thm}

We will prove Theorem \ref{thm:equivalences:1hmin} in Section \ref{sec:equiv:proof}. Until then, we continue distinguishing between the four notions.

\subsection{Basic results under 1-h$^\mix$-minimality}

Many basic results from \cite{CHR} also hold in mixed characteristic with only minor changes. In particular, this includes all the results of Sections~2.4--2.6 from \cite{CHR}, and those of Section~2.8 up to Lemma~2.8.5. In the following, we state the precise version of the mixed characteristic results. When a proof is (almost) identical to the corresponding one from \cite{CHR}, we only give a sketch, so that the reader has the choice between reading the short version here or the long version in \cite{CHR}.

\begin{lem}\label{lem: basics 1hmixmin}
Assume that $\cT$ is $0$-h$^\mix$-minimal.
The following results are true for any model $K$ of $\cT$.
\begin{enumerate}
\item (Adding constants 
\cite[Lemma~2.4.1]{CHR}). If $A\subset   K\cup \RV_{|n|}\eq$ then $\Th_{\cL(A)}(K)$ is $0$-h$^\mix$-minimal. (And similarly, if $\cT$ is $\ell$-h$^\mix$-minimal, then $\Th_{\cL(A)}(K)$ is $\ell$-h$^\mix$-minimal for any $\ell$ and any such $A$.)
\item (Preparation is first order 
\cite[Lemma~2.4.2]{CHR}). If $X_q, C_q$ are $\emptyset$-definable families of subsets of $K$ with $q$ running over an $\emptyset$-definable $Q$ in an arbitrary imaginary sort, and $C_q$ is finite for all $q$, then the set of $(q, \lambda)\in Q\times \Gamma_K^\times$ with $\lambda\leq 1$ such that $C_q$ $\lambda$-prepares $X_q$ is $\emptyset$-definable.
\item ($\exists^\infty$-elimination 
\cite[Lemma~2.5.2]{CHR}). Every infinite $\emptyset$-definable $X\subset   K$ contains an open ball.
\item (Finite sets are $\RV$-parametrized 
\cite[Lemma~2.5.3]{CHR}). If $C_q\subset   K$ is an $\emptyset$-definable family of finite sets, for $q$ in some arbitrary $\emptyset$-definable imaginary sort $Q$, then there is a $\emptyset$-definable family of injections $f_q: C_q\to \RV_{|n|}^k$ (for some $k \ge 0$ and $n \ge 1$).
\item ($\lambda$-next balls as unions of fibres 
\cite[Lemma~2.5.4]{CHR}). Let $C\subset   K$ be a finite $\emptyset$-definable set.  Then for any $\lambda \le 1$, there is a $\{\lambda\}$-definable map $f: K\to \RV_{|n|}^k\times \RV_{|m|\lambda}$ (for some $k$, $n$ and $m$) such that every ball $\lambda$-next to $C$ is a union of fibres of $f$. If moreover, $\lambda = |n'|$ for some integer $n'>0$ then we can ensure that for some integer $p>0$, every ball $|p|$-next to $C$ is contained in a fibre of $f$.
\end{enumerate}
\end{lem}

In (1) of Lemma \ref{lem: basics 1hmixmin}, by $\RV_{|n|}\eq$, we mean imaginary sorts of the form $(\RV_{|n|})^m/\mathord{\sim}$, for some $m$ and some $\emptyset$-definable equivalence relations $\sim$.

\begin{proof}[Proof of Lemma~\ref{lem: basics 1hmixmin}]
(1) and (2) are straightforward from the definition.
(3) also follows directly by preparing $X$.

(4) Using (3), we can assume that $\#C_q$ does not depend on $q$. We define $a_q := \frac1{\#C_q}\sum_{x \in C_q} x$ and $\hat f_q\colon C_q\to \RV_{|n|}, x \mapsto \rv_{|n|}(x - a_q)$ for some $n$ which is a multiple of $\#C_q$. This implies that $\hat f_q$ is not constant on $C_q$. We then apply induction to the family consisting of all fibers of $f_q$, for all $q$.
(We take the final $n$ to be a multiple of all cardinalities $\#C_q$ appearing during this process.)

(5) For $x\in K$, define
\begin{align*}
\mu(x) &= \min\{|x-c|\mid c\in C\},\\
C(x) &= \{c\in C\mid |x-c| = \mu(x)\}, \\
a(x) &= \frac{1}{\#C(x)}\sum_{c\in C(x)}c.
\end{align*}
The map $a: K\to K$ has finite image, so by (4) we can find an injection $\alpha: \operatorname{im} a\to \RV_{|n|}^k$. Let $m = \max\{C(x)\mid x\in K\}!$ and define the map $f$ as
\[
f(x) := \begin{cases}(\alpha(a(x)), \rv_{|m|\lambda}(x-a(x))) & \text{ if } |x-a(x)|\geq \mu(x)\lambda|m| \\
	(\alpha(a(x)), \rv_{|m|\lambda}(0) ) & \text{ if } |x-a(x)|< \mu(x)\lambda|m|. \end{cases}
\]
If $x$ and $x'$ are in the same ball $\lambda$-next to $C$ then $C(x)=C(x')$ and so $a(x) = a(x')$. If $x$ and $x'$ are in the same fibre of $f$ and we are in the first case then the fact that $\rv_{|m|\lambda}(x-a(x)) = \rv_{|m|\lambda}(x'-a(x'))$ implies that $\rv_\lambda(x-c) = \rv_\lambda(x'-c)$ for any $c\in C$. Thus any fibre of $f$ is contained in a ball $\lambda$-next to $C$.

If $\lambda = |n'|$ for some integer $n'$, then one can take $p = n'\cdot m$. Then $f$ will be constant on any ball $|p|$-next to $C$.
\end{proof}

The following proposition states that we can also prepare families, similarly to 
\cite[Proposition~2.6.2]{CHR}. That proposition will also be needed for $0$-h$^\mix$-minimality, so we formulate it more generally.

\begin{prop}[Preparing families]\label{prop: 1hmixmin preparing families}
Assume $\ell$-h$^\mix$-minimality for $\cT$ and some $\ell\ge 0$.
Let $K$ be a model of $\cT$. Let $A\subset  K$ and let
\[
W\subset   K\times \RV^k_{|n|} \times \bigcup_{\lambda \le 1}\RV_\lambda^\ell
\]
be $(A \cup \RV_{|n|})$-definable for some integers $k$ and $n \ge 1$. Then there exists a finite $A$-definable set $C \subset K$ and a positive integer $m$ such that for any $\lambda\leq 1$ and any ball $B$ which is $|m|\lambda$-next to $C$, the set $W_{x, \lambda}:= \{(\xi, \xi') \in \RV^k_{|n|} \times \RV_\lambda^\ell \mid (x, \xi, \xi') \in W\}$ is independent of $x$ as $x$ runs over $B$.
\end{prop}

\begin{proof}
For each $\lambda$ and each $(\xi, \xi') \in \RV^k_{|n|} \times \RV_\lambda^\ell$, let $C_{\lambda,\xi,\xi'}$ be a finite $A$-definable set $|m_{\lambda,\xi,\xi'}|\lambda$-preparing the fiber $W_{\xi,\xi'} \subset K$.
Using compactness, we may assume that for varying $\lambda,\xi,\xi'$, there are only finitely many different sets $C_{\lambda,\xi,\xi'}$ and integers $m_{\lambda,\xi,\xi'}$. Let $C$ be the union of the $C_{\lambda,\xi,\xi'}$ and $m$ be the least common multiple of the $m_{\lambda,\xi,\xi'}$.
\end{proof}

\begin{remark}\label{rem: allow eq}
In that proposition, instead of $\RV^k_{|n|}$, we can have any product $Z$ of sorts from $\RV_{|n|}\eq$ (including in particular $\Gamma_K$). Indeed, given such a $W \subset K \times Z \times \bigcup_{\lambda \le 1}\RV_\lambda^\ell$, we can apply Proposition~\ref{prop: 1hmixmin preparing families} to the preimage of $W$ in $K \times \RV^k_{|n|} \times \bigcup_{\lambda \le 1}\RV_\lambda^\ell$ under some quotient map $\RV^k_{|n|} \to Z$. The same also applies to \Cref{cor: 1hmixmin preparing families} below.
\end{remark}

We will mostly apply the following special case of Proposition~\ref{prop: 1hmixmin preparing families}:

\begin{cor}\label{cor: 1hmixmin preparing families}
Let $K$ be a model of a $0$-h$^\mix$-minimal theory $\cT$, let $A\subset K$ and let
\[
W\subset   K\times \RV^k_{|n|}
\]
be $(A \cup \RV_{|n|})$-definable, for some integers $k$ and $n \ge 1$. Then there exists a finite $A$-definable set $C \subset K$ and a positive integer $m$ such that for any ball $B$ which is $|m|$-next to $C$, the set $W_{x}:= \{\xi \in \RV^k_{|n|}  \mid (x, \xi) \in W\}$ is independent of $x$ as $x$ runs over $B$.
\end{cor}

Note that in the corollary, instead of saying that $W_x$ is constant on each $B$, one could equivalently also say that for each $\xi \in \RV^k_{|n|}$, the set $W_{\xi}:= \{x \in K \mid (x, \xi) \in W\}$ is $|m|$-prepared by $C$. In applications, we will sometimes use this point of view without further notice.

Yet another point of view of the corollary is obtained if $W$ is the graph of a function $f\colon K \to \RV^k_{|n|}$. In that case, the conclusion is that $f$ is constant on each ball $|m|$-next to $C$.



We also obtain the following corollary about preparing families in in all models $K$ of the (possibly non-complete) theory $\cT$. The point here is that the integer $m$ can be taken uniformly over all models.

\begin{cor}[of \Cref{prop: 1hmixmin preparing families}]\label{cor:prep:nc}
Assume that $\cT$ is $\ell$-\hmix-minimal for some $\ell\ge 0$, and suppose that $\phi$ is an $\cL$-formula such that for every model $K \models \cT$, $W_K := \phi(K)$
is a subset of $K\times \RV^k_{|n|}\times \RV_{\lambda_K}^\ell$ for some $k$, $n$, and some $\lambda_K \le 1$ in $\Gamma_K^\times$.
Then there exists an $\cL$-formula $\psi$ and an integer $m \ge 1$ such that for every model $K \models \cT$,
$C_K := \psi(K)$ is a finite subset of $K$ which $\lambda_K\cdot |m|$-prepares $W_K$ in the following sense:
For every ball $B \subset K$ which is $\lambda_K\cdot |m|$-next to $C_K$,
the fiber
$$
W_{K,x} := \{\xi \in \RV^k_{|n|} \times \RV_{\lambda_K}^\ell \mid (x,\xi)\in W_K\}
$$
does not depend on $x$ when $x$ runs over $B$.
\end{cor}

\begin{proof}
Let $\phi$ be given as in the statement. Whether a pair $(m, \psi)$ works as desired in a model $K$ can be expressed by
an $\cL$-sentence, by \Cref{lem: basics 1hmixmin}. By compactness and $\ell$-\hmix-minimality,
we deduce that there exist finitely many pairs $(m_i, \psi_i)$ which cover all models. We may furthermore assume that the sets $\psi_i(K)$ are finite for each model $K$.  We are done by letting  $m$ be the least common multiple of the $m_i$
(so that $|m| \le |m_i|$ for each $i$) and $\psi$ be the disjunction of the $\psi_i$.
\end{proof}

Proposition~\ref{prop: 1hmixmin preparing families} also implies the following corollaries.
Denote by $\RV_\bullet$ the disjoint union of all $\RV_\lambda$ for $\lambda\leq 1$.

\begin{cor}\label{cor: 1hmixmin RVfiniteness}
Assume 1-h$^\mix$-minimality for $\cT$.
The following hold for any model $K$ of $\cT$.
\begin{enumerate}
\item ($\RV$-unions stay finite 
\cite[Corollary~2.6.7]{CHR}). Let $W\subset   K \times \RV_\bullet^{k}$ be $\emptyset$-definable such that $W_\xi$ is finite for any $\xi\in   \RV_\bullet^{k}$. Then the union $\bigcup_\xi W_\xi$ is also finite.
\item (Finite image in $K$ 
\cite[Corollary~2.6.8]{CHR}). The image of any $\emptyset$-definable $f: \RV_\bullet^{k}\to K$ is finite.
\item (Removing $\RV$-parameters 
\cite[Corollary~2.6.10]{CHR}). Let $C$ be a finite $A\cup \RV_{|n|}\eq$-definable set, for some $A\subset  K$ and some integer $n>0$. Then there exists a finite $A$-definable set $C'$ containing $C$.
\end{enumerate}
\end{cor}

\begin{proof}
(1) The case $k = 1$ follows by applying Proposition~\ref{prop: 1hmixmin preparing families}. (Each $W_\xi$ is contained in the set $C$ one obtains.) For $k \ge 2$, use induction (adding parameters to the language using Lemma~\ref{lem: basics 1hmixmin} (1)).

(2) Apply (1) to the graph of $f$.

(3) $C$ is a fiber of an $A$-definable subset $W \subset K \times \RV_\bullet^k$, for some $k$. Apply (1) to $W$.
\end{proof}

To obtain  (T1,T2) from 1-h$^\mix$-minimality for $\cT$, we follow 
\cite[Section~2.8]{CHR}. Here we have to be slightly more careful in our formulations and proofs. First, we set some notation.

\begin{defn}
If $B$ is an open ball in $K$, and $\lambda$ is in $\Gamma_K^\times$, $\lambda\leq 1$, then a \emph{$\lambda$-shrinking of $B$} is an open ball $B'\subset B$ with
\[
\radop B' = \lambda \radop B.
\]
 \end{defn}

\begin{lem}\label{lem: 1hmixmin basic function props}
Assume 1-h$^\mix$-minimality for $\cT$.
Let $K$ be a model of $\cT$ and let $f: K\to K$ be an $\emptyset$-definable function. Then the following hold.
\begin{enumerate}
\item\label{eq.dim.pres} (Basic preservation of dimension 
\cite[Lemma~2.8.1]{CHR}). The set of $y\in K$ for which $f^{-1}(y)$ is infinite, is finite.
\item (Piecewise constant or injective 
\cite[Lemma~2.8.2]{CHR}). There exists a finite $\emptyset$-definable set $C$ and a positive integer $m$ such that on any ball $B$ $|m|$-next to $C$, $f$ is either constant or injective.
\item (Images of most balls are almost balls 
\cite[Lemma~2.8.3]{CHR}). There exists a finite $\emptyset$-definable set $C$ and positive integers $m,n$ such that for any open ball $B$ contained in a ball $|m|$-next to $C$, either $f(B)$ is a singleton, or for any $y\in f(B)$, there are open balls $B', B''$ for which $y\in B'\subset   f(B)\subset   B''$ and
\[
\radop B' \geq |n| \radop B''.
\]
\item (Preservation of scaling factor 
\cite[Lemma~2.8.4]{CHR}). Suppose that there are $\alpha, \beta$ in $\Gamma_K^\times$, $\alpha < 1$, such that for every open ball $B\subset  \cM_K$ of radius $\alpha$, the image $f(B)$ is contained in an open ball of radius $\beta$. Assume moreover that there is an integer $p$ and open balls $B', B''$ such that $B'\subset   f(\cM_K)\subset   B''$ and
\[
\radop B' \geq |p| \radop B''.
\]
Then $f(\cM_K)$ is contained in an open ball of radius at most $\frac{\beta}{|n|\alpha}$ for some positive integer $n$.
\end{enumerate}
\end{lem}
\begin{proof}
(1) Suppose for contradiction that $f$ has infinitely many infinite fibers. Let $X \subset K$ be the subset of the domain of $f$ where $f$ is locally constant,
i.e, of points $x \in K$ such that $f$ is constant on $B_{<\lambda}(x)$ for some $\lambda \in \Gamma_K^\times$.
Since each infinite fiber of $f$ contains a ball (by preparation of the fiber), the restriction of $f$ to $X$ still has infinite image.

Let $C \subset K$ be a finite set $|m|$-preparing $X$, for some integer $m \ge 1$. Then enlarge $C$ and $m$ in such a way that for each $\lambda \in \Gamma_K^\times$, whether $f$ is constant on $B_{<\lambda}(x)$ only depends on the ball $|m|$-next to $C$ containing $x$. This is possible by applying Corollary~\ref{cor: 1hmixmin preparing families} (and Remark~\ref{rem: allow eq})
to the set $W \subset K \times \Gamma_K^\times$ of those $(x,\lambda)$ for which $f$ is constant on $B_{<\lambda}(x)$.

By Corollary~\ref{cor: 1hmixmin RVfiniteness} (1), there exists a ball $B_0 \subset X$ $|m|$-next to $C$ such that $f(B_0)$ is still infinite. Fix such a $B_0$, and also fix a $\lambda_0 \in \Gamma_K^\times$ such that $f$ is constant on every open ball of radius $\lambda_0$ contained in $B_0$.

The family of all those balls can definably be parametrized by a subset $Z \subset \RV_{\lambda_1}$, for $\lambda_1 = \lambda_0/\radop B_0$. Now $f$ induces a map from $Z$ to $K$ with infinite image, contradicting Corollary~\ref{cor: 1hmixmin RVfiniteness} (1).

\medskip

(2) Firstly, find $C$ and $m$ such that $C$ $|m|$-prepares every infinite fiber of $f$ (using (1) to see that there are only finitely many infinite fibers). Then apply Lemma~\ref{lem: basics 1hmixmin} (4) to the family of finite fibers of $f$, to obtain an injective map from each finite fiber of $f$ to $\RV^k_{|n|}$ for some $n$ and $k$.
We put all those maps together to one single map
$g\colon Y \to \RV^k_{|n|}$, where $Y$ is
the union of all finite fibers of $f$. Enlarge $C$ and $m$
in such a way that $g$ is constant on each ball $B \subset Y$ $|m|$-next to $C$ (by applying Corollary~\ref{cor: 1hmixmin preparing families} to the graph of $g$). Then, for each ball $B$ $|m|$-next to $C$, either $B$ is entirely contained in an infinite fiber of $f$ (and hence $f$ is constant on $B$, as desired), or $B \subset Y$ and $g$ is constant on $B$. In that case, $f$ is injective on $B$, since $f(x_1) = f(x_2)$ for $x_1, x_2 \in B$ would imply that $x_1$ and $x_2$ lie in the same fiber of $f$, contradicting that $g$ is injective on each fiber of $f$.

\medskip

(3) Use (2) to obtain a finite $\emptyset$-definable set $C$ and an integer $m$ such that on any ball $|m|$-next to $C$, $f$ is either constant or injective. Let $W_0\subset   K\times (\Gamma_K^\times)^2$ consist of those $(x, \lambda, \mu)$ with $\mu\leq 1$ such that
for every $y\in f(B_{<\lambda}(x))$ there are open balls $B', B''$ with $y\in B'\subset   f(B_{<\lambda}(x))\subset   B''$ and
\[
\radop B' \geq \mu \radop B''.
\]
We enlarge $C$ and $m$, such that $C$ also $|m|$-prepares this set $W_0$.

The set $C$ is already as desired, but $m$ will later be enlarged to some $m\cdot m''$. Note that we already simplified the statement we will need to prove, namely: Firstly, it suffices to consider balls $B$ which are contained in a ball $B_1$ $|m|$-next to $C$ on which $f$ is injective; and secondly, for each $\lambda \le 1$, it suffices to find a single $\lambda\cdot|m''|$-shrinking $B$ of $B_1$ for which the lemma holds (using some $n$ which we still need to specify). Indeed, the fact that $C$ $m$-prepares $W_0$ implies that then, the lemma also holds for all translates of $B$ within $B_1$ (using the same $n$).

Before we can continue, we need to do some preparation on the range side of $f$: We want to find a finite $\emptyset$-definable $D$ and a positive integer $p$ such that for every $\lambda \le 1$ and every ball $B$ $\lambda|m|$-next to $C$, the set $f(B)$ is $\lambda|p|$-prepared by $D$.
To see that such a $D$ exists, first note that
by Lemma~\ref{lem: basics 1hmixmin} (5), there
exists a $\lambda$-definable map $g_\lambda\colon K \to \RV_{|n'|}^k\times \RV_{|m'|\cdot \lambda}$ such that every ball $\lambda|m|$-next to $C$ is a union of fibers of $g_\lambda$. We may assume (by compactness) that the maps $g_\lambda$ form a $\emptyset$-definable family.
Now apply Proposition \ref{prop: 1hmixmin preparing families} to the set
\[
W := \{(f(x), \lambda, g_\lambda(x)) \in K \times \Gamma_K^\times \times \RV_{|n'|}^k \times  \bigcup_\mu \RV_{\mu} \mid x \in K, \lambda \le 1\}
\]
and let $D$ and $p$ be the result.
To see that this works, let $\lambda\le 1$ and $B$ $\lambda|m|$-next do $C$ be given. We have $B = g_\lambda^{-1}(\Xi)$ for $\Xi := g_\lambda(B)$,
and the image $f(B)$ consists of those $y \in K$ for which the fiber $W_{(y, \lambda)}$ is not disjoint from $\Xi$. Since this fiber is constant when $y$ runs over a ball $B'$ $\lambda|p|$-next to $D$, we either have $B' \subset f(B)$ or $B' \cap f(B) = \emptyset$, as desired.

Now that we constructed $D$ and $p$ (on the range side of $f$), we construct another set $C'$ and integer $m'$ on the domain side of $f$:
For $\lambda\in \Gamma_K^\times$, $\lambda\leq 1$, use \Cref{lem: basics 1hmixmin} to get a map $h_\lambda: K\to \RV_{|p|}^k\times \RV_{\lambda|n''|}$ such that each fibre is contained in a ball $\lambda|p|$-next to $D$. We may assume that $h_\lambda$ is an $\emptyset$-definable family of maps, with a parameter $\lambda$. Let $C'$ be a finite $\emptyset$-definable set which $\lambda|m'|$-prepares every map $h_\lambda\circ f$. This can again be done using Proposition \ref{prop: 1hmixmin preparing families}.

We need one last ingredient before we can verify that the lemma holds: Since $C'$ is finite, there exists an integer $p' \ge 1$ such that every ball $B_1$ $|m m'|$-next to $C$ has a $|p'|$-shrinking $B_1' \subset B_1$ disjoint from $C'$.

We now claim that the lemma holds using $C$, $|mm'p'|$, and $n = m'p'$. According to the beginning of the proof, it suffices to check that given a ball $B_1$ $|m|$-next to $C$ on which $f$ is injective, and given a $\lambda \le 1$,
there exists a $\lambda \cdot |p'  m'|$-shrinking $B$ of $B_1$
for which the claim holds.
Choose a $|p'|$-shrinking $B_1'$ of $B_1$ disjoint from $C'$ and choose $B$ to be any $\lambda|m'|$-shrinking of $B_1'$. To finish the proof, we show that the lemma holds for this $B$.

On the one hand, since $B$ is a $\lambda \cdot|m'|$-shrinking of $B_1'$ (and $B_1'$ is disjoint from $C'$), $B$ is contained in a ball $\lambda \cdot |m'|$-next to $C'$.
By definition of $C'$ this means that $h_{\lambda} \circ f$ is constant on $B$ and hence (by definition of $h_{\lambda}$) $f(B)$ is contained in a ball $B''$ $\lambda|p|$-next to $D$.

On the other hand, $B$ is $\lambda|mm'p'|$-next to $C$, so that $f(B)$ is $\lambda|pm'p'|$-prepared by $D$ (by definition of $D$). Thus for any $y\in f(B)$, we obtain that the entire ball $B'$ $\lambda|pm'p'|$-next to $D$ containing $y$ is contained in $f(B)$. This ball $B'$ is just the (unique) $|m'p'|$-shrinking of $B''$ containing $y$; in particular, $\radop B' = |m'p'|\radop B''$, as desired.

\medskip

(4) If the residue field of $K$ has characteristic zero, then this is exactly~\cite[Lemma~2.8.4]{CHR}. So we may assume that we are in mixed characteristic (though the following proof easily also adapts to the equicharacteristic $0$ case.) The family of radius $\alpha$ open balls in $\cM_K$ can be definable parametrized by the set $\Lambda = \rv_\alpha(\cM_K)$. Namely, for $\xi\in \Lambda$ let $B_\xi\subset \cM_K$ be the unique open ball of radius $\alpha$ containing $\rv_\alpha^{-1}(\xi)$. Using 1-h$^\mix$-minimality, we can find a finite set $C\subset K$ and a positive integer $m$ such that $C$ $|m|\alpha$-prepares every set $f(B_\xi)$, for every $\xi\in \Lambda$. In other words, if $\xi\in \Lambda$, and if $B'$ is an open ball $|m|\alpha$-next to $C$ then either $B'\subset f(B_\xi)$ or $B'\cap f(B_\xi) = \emptyset$. Now, if $B'$ has radius strictly larger than $\beta$, then it follows that $B'\cap f(B_\xi) = \emptyset$. Therefore, $f(\cM_K)$ is contained in the union of $C$ and all balls $|m|\alpha$-next to $C$ of radius at most $\beta$. This union is equal to the finite union of all closed balls of the form $B_{\leq \beta/(|m|\alpha)}(c), c\in C$. In particular the open ball $B'$ is contained in this finite union of closed balls of radius $\beta/(|m|\alpha)$. Since we are in mixed characteristic, there exists some positive integer $q$ such that $\radop B'\leq \beta/(|q|\alpha)$. But then $B''$ is an open ball containing $f(\cM_K)$ of radius at most $\beta / (|pq|\alpha)$, finishing the proof. 
\end{proof}

Finally, we can prove an approximate valuative Jacobian property in mixed characteristic. The lemma and its proof are similar to 
\cite[Lemma~2.8.5]{CHR}
and is sharpened to an actual Jacobian property below in Corollary \ref{cor:T1exactmix}. Note that the sharpened version is obtained only using a huge detour, involving approximations by second degree Taylor polynomials. We do not see a more direct proof of this sharpened version.

\begin{lem}[{Approximate valuative Jacobian property}] \label{lem: 1hmix valjac}
Assume that $\cT$ is 1-h$^\mix$-minimal.
Let $K$ be a model of $\cT$ and let $f: K\to K$ be an $A$-definable function, for some $A\subset K\cup \RV_{|n|}$. Then there exists a finite $A$-definable set $C$ and a positive integer $m$ such that for every ball $B$ $|m|$-next to $C$ there exists a $\mu_B\in \Gamma_K$ such that for all $x, y\in B$ we have
\[
\mu_B|m||x-y| \leq |f(x)-f(y)| \leq \mu_B \left|  \frac{1}{m} \right| |x-y|.
\]
\end{lem}
\begin{proof}
We may assume that $A=\emptyset$ by Lemma \ref{lem: basics 1hmixmin}. Using Lemma \ref{lem: 1hmixmin basic function props} (3), we can find a finite $\emptyset$-definable set $C_0$ and positive integers $m,n$ such that
\begin{itemize}
\item $f$ is constant or injective on balls $|m|$-next to $C_0$, and
\item if $B$ is an open ball contained in a ball $|m|$-next to $C_0$ then either $f(B)$ is a singleton, or there are open balls $B'\subset   f(B)\subset   B''$ such that the radii of $B'$ and $B''$ differ by at most $|n|$.
\end{itemize}
For an open ball $B = B_{<\alpha}(x)$ contained in a ball $|m|$-next to $C$ on which $f$ is injective, define $\mu(x, \alpha)$ to be the (convex) set of $\mu\in \Gamma_K^\times$ for which $f(B)$ is contained in an open ball of radius $\mu$ and contains an open ball of radius $|n|\mu$. Note that we have $|n|\mu(x,\alpha)\leq \mu(x,\alpha)$, in the sense that
for every $\nu \in |n|\mu(x,\alpha)$ and every $\nu' \in \mu(x,\alpha)$, we have $\nu \leq \nu'$. (In the following, inequalities between convex subsets of $\Gamma_K$ are always meant in this sense.) Also define
\[
s(x,\alpha) = \{\mu/\alpha \mid \mu \in \mu(x,\alpha)\}.
\]
(An element $\nu \in s(x,\alpha)$ is a kind of ``scaling factor'': the ball $B_{<\alpha}(x)$ is sent into a ball of radius $\alpha\nu$.)

We will now enlarge $C$ and $m$ to obtain the following claim. As a side remark, note that we will keep $n$ fixed for the entire proof, so that the definitions of $\mu(x,\alpha)$ and $s(x, \alpha)$ are unaltered.

\medskip

\paragraph{Claim 1:} By possibly enlarging $C$ and $m$, we can achieve that
\begin{enumerate}
\item $\mu(x,\alpha)$ and $s(x,\alpha)$ are independent of $x$ as $x$ runs over a ball $|m|$-next to $C$,
\item if $B_{<\alpha}(x)$ and $B_{<\beta}(y)$ are open balls contained in the same ball $|m|$-next to $C$ and $\alpha\leq \beta$, then
\[
s(y,\beta) \leq s(x,\alpha)/|m|.
\]
\end{enumerate}

\medskip

To prove Claim~1, let $W\subset   K\times (\Gamma_K^\times)^2$ consist of those $(x, \lambda, \mu)$ such that 
$f(B_{<\lambda}(x))$ is contained in a ball $B'$ of radius $\mu$, and contains a ball of radius $|n|\mu$.
Enlarge $C$ and $m$ such that $C$ $|m|$-prepares this set $W$. Note that item (1) of the claim then holds by preparation. Let $B$ be $|m|$-next to $C$. If $f$ is constant on $B$ then there is nothing to check, so we can assume that $f$ is injective on $B$. By item (1), $\mu(x,\alpha)$ and $s(x,\alpha)$ are constant when $x$ runs over $B$ and $\alpha\leq \radop B$ is fixed, so we write simply $\mu(\alpha)$ and $s(\alpha)$. It remains to prove item (2) of Claim~1 (after possibly enlarging $m$ once more). Fix $\mu\in \mu(\alpha)$ and take $\alpha\leq \beta\leq \radop B$. Then any ball of radius $\alpha$ inside $B$ has image under $f$ contained in a ball of radius at most $\mu$. Hence, by a rescaled version of Lemma \ref{lem: 1hmixmin basic function props}(4), there exists an integer $p \ge 1$ such that
for $x \in B$, the image $f(B_{<\beta}(x))$ is contained in
an open ball of radius $\frac{\mu\cdot \beta}{|p| \cdot\alpha}$. In particular, we have
\[
|n|\mu(\beta)\leq \frac{\mu(\alpha)\cdot \beta}{|p| \cdot\alpha}.
\]
But this means precisely that
\[
s(\beta)\leq s(\alpha)/|np|,
\]
which proves Claim~1 (after replacing $m$ by e.g.\ $mnp$). 

\medskip

To prove the lemma, we also need an inequality opposite to the one of Claim~1 (2). More precisely:

\medskip

\paragraph{Claim 2:} After possibly further enlarging $C$ and $m$, in Claim~1 (2), we can additionally obtain
\[
s(y,\alpha) \leq s(x,\beta)/|m|.
\]

\medskip

The idea of the proof of Claim~2 is to apply Claim~1 to $f^{-1}$ (which in reality only exists piecewise). This is made precise as follows:

Denote by $Y$ the set of $y\in K$ for which $f^{-1}(y)$ is finite. This is a cofinite $\emptyset$-definable set in $K$, by Lemmas
\ref{lem: basics 1hmixmin} and \ref{lem: 1hmixmin basic function props} (1). Use Lemma \ref{lem: basics 1hmixmin} (4) to obtain a $\emptyset$-definable family of injections
\[
h_y: f^{-1}(y)\to \RV_{|n_0|}^k,
\]
for $y\in Y$. For $\eta\in \RV_{|n_0|}^k$ define
\begin{align*}
Y_\eta &= \{y\in Y\mid \eta\in \operatorname{im} h_y\}, \\
g_\eta &: Y_\eta \to K: y\mapsto h_y^{-1}(\eta).
\end{align*}
Then we have that
\[
\bigcup_{y\in Y}\{y\}\times f^{-1}(y) = \bigsqcup_{\eta\in \RV_{|n_0|}^k}\operatorname{graph}(g_\eta).
\]
For each $\eta$, apply \Cref{lem: 1hmixmin basic function props} (3) to $g_\eta$ (extended by $0$ outside of $Y_\eta$) to obtain a finite $\eta$-definable set $D_\eta$ and integers $m_\eta, n_\eta$.
By compactness, we may take $n' := n_\eta$ independent of $\eta$. Now enlarge $D_\eta$ and $m_\eta$
using Claim~1, so that (1) and (2) hold for $g_\eta$ (where the corresponding $\mu$ and $s$ are defined using $n'$). We may moreover assume that $D_\eta$ $|m_\eta|$-prepares $Y_\eta$ (enlarging $D_\eta$ and $m_\eta$ once more, if necessary).
After that, apply compactness once more to make
$m':=m_\eta$ independent of $\eta$ and to turn $(D_\eta)_\eta$ into a $\emptyset$-definable family. By Corollary \ref{cor: 1hmixmin RVfiniteness} the union $D=\bigcup_\eta D_\eta$ is a finite $\emptyset$-definable set.

Using Lemma \ref{lem: basics 1hmixmin} (5), take a $\emptyset$-definable function $\chi: K\to \RV_{|n_1|}^{k'}$ such that every ball $|m'|$-next to $D$ is a union of fibres of $\chi$.
Then choose $C_2$ and $C_3$ such that, after possibly enlarging $m$, $C_2$ $|m|$-prepares the family of sets $f^{-1}(\chi^{-1}(\eta))$ for $\eta\in \RV_{|n_1|}^{k'}$, and $C_3$ $|m|$-prepares the family of images $g_\eta(Y_\eta)$.
We will now prove that Claim~2 holds after replacing $C$ by $C \cup C_2 \cup C_3$ (and some further enlargement of $m$).

So suppose that we have open balls $B_1, B$ with $\radop B_1 = \alpha \leq \radop B = \beta$ which are contained in the same ball $|m|$-next to $C$. If $f$ is constant on $B$ then we are done, so assume that $f$ is injective on $B$. We may assume that $B_1\subset B$ since $\mu(x, \alpha)$ is independent of $x$ as $x$ runs over $B$. Let $B'\subset f(B)\subset B''$ be open balls whose radii differ by at most $|n|$. By definition of $C_2$, $\chi \circ f$ is constant on $B$, so (by definition of $\chi$), $f(B)$ is contained in a ball $|m'|$-next to $D$. Perhaps after shrinking $B''$, we can assume that also $B''$ is contained in a ball $|m'|$-next to $D$. By definition of $C_3$, there is a (unique) $\eta\in \RV_{|n_0|}^k$ such that $B \subset g_\eta(Y_\eta)$; this implies that $f|_B$
and $g_\eta|_{f(B)}$ are mutually inverse bijections between $B$ and $f(B)$.
Using that we applied \Cref{lem: 1hmixmin basic function props} to $g_\eta$, take open balls $\tilde{B}'\subset g_\eta(B')\subset \overline{B}'$ whose radii differ by at most $|n'|$ and do the same with open balls $\tilde{B}''\subset g_\eta(B'')\subset \overline{B}''$.
Note that we have a chain of inclusions
\[
\tilde{B}' \subset g_\eta(B') \subset g_\eta(f(B)) = B \subset g_\eta(B'') \subset  \overline{B}''.
\]
Choose similar balls corresponding to $B_1$: $B_1', B_1'', \tilde{B}_1', \dots$. We may certainly assume that $\radop B_1' \le \radop B''$. Therefore,
our application of Claim~1 to $g_\eta$ yields that
\[
\frac{\radop \overline{B}''}{\radop B''}
\le
\frac{\radop \overline{B}'_1}{|m'|\radop B'_1}.
\]
Combining this with $\radop B \le \radop\overline{B}''$,
$|n'|\radop \overline{B}'_1 \le \radop\tilde{B}'_1 \le \radop B_1$ and $|n|\radop B''_1 \le \radop B'_1$,
we obtain
\[
\frac{\radop B}{\radop B''}
\le
\frac{\radop \overline{B}''}{\radop B''}
\le
\frac{\radop \overline{B}'_1}{|m'|\radop B'_1}
\le
\frac{\radop B_1}{|n'm'|\radop B'_1}
\le
\frac{\radop B_1}{|nn'm'|\radop B''_1}.
\]
Since
$\frac{\radop B''_1}{\radop B_1} \in s(\alpha)$ and
$\frac{\radop B''}{\radop B} \in s(\beta)$, we deduce
\[
s(\alpha) \le s(\beta)/|n^3n'm'|
\]
(where the additional factor $n^2$ takes into account the length of the convex sets $s(\alpha)$ and $s(\beta)$).
This finishes the proof of Claim~2 (where we take $m$ to be a multiple of $n^3n'm'$).

\medskip

We are now ready to prove the lemma itself.
Take $x,y$ in the same ball $B$ $|m|$-next to $C$. Denote by $\beta$ the open radius of $B$. Let $\mu\in \mu(x, \beta)$, so that $f(B)$ is contained in a ball of radius $\mu$ and contains a ball of radius $|n|\mu$.
We will show that we can take $\mu_B = \mu/\beta$.

If we choose $\alpha > |x-y|$, then $x$ and $y$ are contained in an open ball of radius $\alpha$. If we moreover choose $\mu'\in \mu(x,\alpha)$, then $f(x), f(y)$ are contained in an open ball of radius $\mu'$. Thus
\[
|f(x)-f(y)|/\alpha < \mu'/\alpha \leq \frac{1}{|m|} \mu /\beta.
\]
Since this holds for any $\alpha > |x-y|$, this gives that
\[
|f(x)-f(y)| \leq \frac{1}{|m|} \frac{\mu}{\beta} |x-y|.
\]
For the other inequality, set $\alpha = |x-y|$ and denote by $B'$ the open ball of radius $\alpha$ around $x$. By the injectivity of $f$ on $B$, $f(y)$ is not in $f(B')$. But $f(B')$ contains an open ball of radius $|n|\mu'$ (for $\mu' \in \mu(x, \alpha)$, as before) and so
\[
|f(x)-f(y)|/\alpha \geq |n|\mu'/\alpha \geq |nm|\mu/\beta.
\]
Thus
\[
|f(x)-f(y)| \geq |mn|\frac{\mu}{\beta}|x-y|.
\]
This proves the Lemma.
\end{proof}

We now obtain already a part of Theorem \ref{thm:equivalences:1hmin}.

\begin{cor}\label{cor.T1.T2}
Any 1-h$^\mix$-minimal theory $\cT$ satisfies (T1,T2).
\end{cor}
\begin{proof}
This is direct from \Cref{eq.dim.pres} from \Cref{lem: 1hmixmin basic function props} and \Cref{lem: 1hmix valjac}.
\end{proof}

\subsection{Basic results under (T1,T2)} 

We now prove some consequences of $\cT$ satisfying (T1,T2) in the sense of Definition~\ref{defn:T1T2}. First note that by applying (T1$^\mix$) to the characteristic function of an $A$-definable set $X \subset K$, for some $A \subset K \cup \RV_{|n|}$, we obtain a finite $A$-definable set $C \subset K$ $|m|$-preparing $X$ for some $m\ge 1$. In Lemma~\ref{lem: property 3 basic props}, we will see that in
(T1$^\mix$), we can even take $C$ to be $(A \cap K)$-definable, so that $\cT$ is $0^\mix$-h-minimal.
For the moment, note that the above weaker statement already suffices to obtain $\exists^\infty$-elimination (namely, (3) of \Cref{lem: basics 1hmixmin}).

\begin{lem}\label{lem: property 3 basic props}
Assume that $\cT$ satisfies (T1,T2).
Let $K$ be a model of $\cT$.
Then
\begin{enumerate}
\item ($\RV_\lambda$-unions stay finite 
\cite[Lemma~2.9.4]{CHR}). If $C_\xi\subset   K$ is a definable family (with parameters) of finite sets, parametrized by $\xi\in \RV_\lambda^k$, then $\bigcup_\xi C_\xi$ is still finite.
\item (Eliminating $\RV$-parameters 
\cite[Lemma~2.9.5]{CHR}). If $f: K\to K$ is $A$-definable for some $A\subset   K\cup \RV_{|n|}$ then we can find a finite $(A\cap K)$-definable set $C$ and an integer $m$ such that $f$ satisfies property (T1$^\mix$) with respect to $C$ and $m$. In particular, $\cT$ is 0-\hmix-minimal.
\end{enumerate}
\end{lem}

The following proof is the same as in \cite{CHR}.

\begin{proof}[Proof of Lemma~\ref{lem: property 3 basic props}]
(1) We induct on $k$, so first assume that $k=1$. Using $\exists^\infty$-elimination, we may assume that the cardinality of $C_\xi$ is constant, say, equal to $m$. Let $\sigma_1, \dots, \sigma_m$ be the elementary symmetric polynomials in $m$ variables, considered as functions on $m$-element-subsets of $K$. Then the map $K^k \to K, x \mapsto \sigma_i(C_{\rv_\lambda(x)})$ is locally constant everywhere except possibly at $0$, so by (T2), it has finite image. Since $\sigma_1(C_\xi), \dots, \sigma_m(C_\xi)$ together determine $C_\xi$, there are only finitely many different $C_\xi$.

For arbitrary $k$, consider $C_{\xi}$ as a definable family $C_{\xi_1, \xi_2}$ with $\xi_1\in \RV_\lambda$ and $\xi_2\in \RV_\lambda^{k-1}$. Then by induction on $k$ the union $\cup_{\xi} C_{\xi} =  \cup_{\xi_1\in \RV_\lambda}\cup_{\xi_2\in \RV_\lambda^k} C_{\xi_1, \xi_2}$ is finite.

(2) Using (T1$^\mix$), we find an $A$-definable $C$. We consider $C$ as an $(A\cap K)$-definable family $C_\xi$, parametrized by $\xi\in \RV_{|n|}^k$. Using (1), the union $C' = \cup_{\xi\in \RV_{|n|}^k} C_\xi$ is $(A\cap K)$-definable, finite, and $f$ satisfies (T1$^\mix$) with respect to $C'$.
\end{proof}

Now that we know that any theory $\cT$ satisfying (T1,T2) is $0^\mix$-h-minimal, we can use \Cref{cor: 1hmixmin preparing families} and Lemma \ref{lem: basics 1hmixmin}.

The following is a first adaptation of 
\cite[Lemma~2.9.6]{CHR}. In \Cref{cor:T1exactmix}, we will obtain an adaptation which is better in the sense that it has a more precise conclusion (and more similar to the equi-characteristic $0$ case).

\begin{lem}[Images of balls]\label{lem:property (3) balls map to balls}
Assume that $\cT$ satisfies (T1,T2). Let $f: K\to K$ be $A$-definable, for some $A\subset   K\cup \RV_{|n|}$. Then there exists a finite $A$-definable $C$ and a positive integer $m$ such that (T$1^\mix$) holds for $f$ and $C$ and $m$ and such that the following holds for any ball $B$ $|m|$-next to $C$. Let $B'\subset   B$ be an open ball and let $\mu_B$ be as in (T$1^\mix$). Then for every $x\in B'$ there are open balls $B_1, \tilde{B}_1$ such that $f(x)\in \tilde{B}_1\subset   f(B')\subset   B_1$ and moreover
\[
|m|\mu_B\radop(B')\leq \radop(\tilde{B}_1)\leq \radop(B_1) \leq \left|\frac{1}{m}\right| \mu_B\radop(B').
\]
\end{lem}

\begin{proof}
Using Lemma \ref{lem: basics 1hmixmin} we may as well assume that $A=\emptyset$. Take a finite $\emptyset$-definable set $C$ and an integer $m$ such that (T$1^\mix$) holds for $f$ with respect to $C$ and $m$. Let $\chi: K\to \RV_{|p|}^k$ be a $\emptyset$-definable map coming from Lemma \ref{lem: basics 1hmixmin}(5) for $C,m$. So every ball $|m|$-next to $C$ is a union of fibres of $\chi$ and every ball $|n'|$-next to $C$ is contained in a fibre of $\chi$. Let $D$ be a finite $\emptyset$-definable set $|n|$-preparing the family $(f(\chi^{-1}(\eta)))_{\eta\in \RV_{|p|}^k}$, here we use \Cref{cor: 1hmixmin preparing families}. Let $\psi: K\to \RV_{|p'|}^{k'}$ be a $\emptyset$-definable map such that every ball $|n|$-next to $D$ is a union of fibres of $\psi$ and every ball $|n''|$-next to $D$ is contained in a fibre of $\psi$. Finally, use \Cref{cor: 1hmixmin preparing families} again to $|q|$-prepare the family $(f^{-1}(\psi^{-1}(\eta)))_{\eta\in \RV_{|p'|}^{k'}}$ with a finite $\emptyset$-definable set $C_0$.

We claim that $C'=C\cup C_0$ suffices. So let $B'$ be a ball $|q|$-next to $C'$ and let $B$ be the ball $|m|$-next to $C$ containing $B'$. We can assume that $\mu_B\neq 0$, for else $f$ is constant on $B$. Fix any open ball $B''$ in $B'$ and let $x,y\in B''$. Then
\[
|f(x)-f(y)|\leq \frac{\mu_B}{|m|}|x-y| < \frac{\mu_B}{|m|}\radop B''.
\]
Therefore, if we denote by $B_1''$ the open ball of radius $\frac{\mu_B}{|m|}\radop B''$ around $f(x)$ then $f(B'')$ is contained in $B_1''$. By definition of $C_0$, $f(B')$ is contained in a ball $B_1$ $|n|$-next to $D$. On the other hand, $f(B)$ is $|n|$-prepared by $D$. Since $B_1\cap f(B)\neq \emptyset$ we have
\[
f(B'')\subset f(B')\subset B_1\subset f(B).
\]
Using property (T$1^\mix$) we see that
\[
\radop B_1\geq |m|\mu_B\radop B'.
\]
Now take $z\in K$ such that $|f(x)-z| < |m|\mu_B\radop B''$. Then we have $z\in B_1\subset f(B)$ and so there is some $x'\in B$ with $f(x')=z$. Applying (T$1^\mix$) one more time yields that
\[
|x-x'|\leq \frac{1}{|m|\mu_B} |f(x)-f(x')| < \radop B''.
\]
We conclude that $f(B'')$ contains the open ball of radius $|m|\mu_B\radop B''$ around $f(x)$.
\end{proof}

\subsection{{Basic results under 1-\heqc-minimality}}\label{ssec:coarsenings.basic}

We provide the tools necessary to transfer preparation results to mixed characteristic, starting from $1$-\heqc-minimality. We prove part of Theorem~\ref{thm:equivalences:1hmin} which states that 1-\heqc-minimality implies 1-\hmix-minimality. 

\begin{notn}
The notion of balls $\lambda$-next to a finite set $C \subset K$
now has different meanings for $|\cdot|$ and for $|\cdot|_\eqc$, with notation from \Cref{not.eqc}. To make clear which of the valuations we mean, we either write $|1|$-next or $|1|_\eqc$-next (instead of just $1$-next).
\end{notn}

\begin{remark}\label{rem:eqc:m}
Suppose that $|\cdot|_{\eqc}$ is non-trivial on $K$. For any $x, x' \in K$, we have
\begin{equation}\label{eq:coarse}
|x|_\eqc \le |x'|_\eqc \quad \iff \quad \exists m \in \NN_{\ge 1}\colon  |m\cdot x| \le |x'|,
\end{equation}
and given a finite set $C \subset K$, the points $x$ and $x'$ lie in the same ball $|1|_\eqc$-next to $C$ if and only if for every integer $m \ge 1$, they lie in the same ball $|m|$-next to $C$.
\end{remark}

Assuming $1$-\heqc-minimality of $K$ as an $\cL_\eqc$-structure, we will be able to find finite $\cL_\eqc$-definable sets. To get back to the smaller language $\cL$, we will use the following lemma:

\begin{lem}[From $\cL_\eqc$-definable to $\cL$-definable]\label{lem:mix}
Let $\cT$ be $1$-\heqc-minimal. Let $K$ be a model of $\cT$ which is $\aleph_0$-saturated and strongly $\aleph_0$-homogeneous as an $\cL\eq$-structure. (Note that this in particular implies that $|\cdot|_\eqc$ is non-trivial.) Then, any finite $\cL_\eqc$-definable set $C \subset K$ is already $\cL$-definable.
\end{lem}

\begin{remark}\label{rem:homog}
It is a standard result that every model $K$ has an elementary extension satisfying the properties of the lemma: Indeed, any model which is special in the sense of \cite[Definition~6.1.1]{TentZiegler} is strongly $\aleph_0$-homogeneous by
\cite[Theorem~6.1.6]{TentZiegler}, and it is easy to construct $\aleph_0$-saturated special models.
\end{remark}

\begin{proof}[Proof of Lemma~\ref{lem:mix}]
It suffices to prove that for any $a \in C$, all realizations of $p:=\tp_{\cL}(a/\emptyset)$ lie in $C$; indeed, by $\aleph_0$-saturation, this then implies that $p$ is algebraic (using that $C$ is finite), and hence isolated by some formula $\phi_p(x)$, and hence $C$ is defined by the disjunction of finitely many such $\phi_p(x)$.

So now suppose for contradiction that there exist $a \in C$ and $a' \in K \setminus C$ which have the same $\cL$-type over $\emptyset$. Then by our homogeneity assumption, there exists an $\cL$-automorphism of $K$ sending $a$ to $a'$ (and hence not fixing $C$ setwise). Such an automorphism also preserves $\cO_{\eqc}$ and hence is an $\cL_\eqc$-automorphism, but this contradicts $C$ being $\cL_\eqc$-definable.
\end{proof}

We obtain yet another part of \Cref{thm:equivalences:1hmin}.

\begin{cor}\label{cor:prep:mix}
Assume that $\cT$ is $1$-\heqc-minimal. Then $\cT$ is 1-h$^\mix$-minimal.
\end{cor}

\begin{proof}
By Remark~\ref{rem:homog}, we may consider a model $K$ of $\cT$ which is sufficiently saturated and sufficiently homogeneous  (as in Lemma~\ref{lem:mix}). Consider an integer $n$, $\lambda\in \Gamma_K^\times$, $A\subset K$,  $\xi \in \RV_\lambda$, $\eta\in \RV_{|n|}^k$ for some integers $k,n>0$, and an $\cL(A\cup \{\eta,\xi\})$-definable set $X\subset K$. We have to show that $X$ can be $|m|\cdot \lambda$-prepared by some finite $\cL(A)$-definable set $C$ for some integer $m>0$.

Let $\lambda_\eqc$ be the image of $\lambda$ in $\Gamma_{K,\eqc}$. Since $B_{<\lambda_\eqc}(1) \subset B_{<\lambda}(1)$, we have a canonical surjection $\RV_{\lambda_\eqc} \to \RV_{\lambda}$. Similarly, there is a canonical surjection $\RV_{\eqc}\to \RV_{|n|}$.
We fix any preimage
$(\xi_\eqc, \eta_\eqc) \in \RV_{\lambda_\eqc} \times \RV_\eqc^k$ of $(\xi,\eta)$, so that $X$ is $\cL(A \cup \{\xi_\eqc, \eta_\eqc\})$-definable.
By $1$-h-minimality for the $\cL_\eqc$-structure on $K$,
there exists a finite $\cL_\eqc(A)$-definable set $C$ such that for every pair $x, x'$ in the same ball $\lambda_\eqc$-next to $C$, we have
$x \in X \iff x' \in X$. By Lemma~\ref{lem:mix}, $C$ is already $\cL(A)$-definable; we claim that it is as desired.
Suppose for contradiction that there exists no $m$ as in the corollary, i.e., for every integer $m \ge 1$, there exists a pair of points $(x, x') \in K^2$ which lie in the same ball $\lambda\cdot|m|$-next to $C$ such that
$x \in X$ but $x' \notin X$.
By $\aleph_0$-saturation (in the language $\cL$), we find a single pair $(x, x') \in K^2$ of points with
$x \in X$ but $x' \notin X$
and which lie in the same ball $\lambda\cdot |m|$-next to $C$ for every $m \ge 1$. The latter implies that $x$ and $x'$
lie in the same ball $\lambda_\eqc$-next to $C$ (by Remark~\ref{rem:eqc:m}), so we get a contradiction to our choice of $C$.
\end{proof}

\subsection{Resplendency}\label{sec:respl}

In order to prove that (3) implies (1) in Theorem \ref{thm:equivalences:1hmin}, we will need a way to add the coarsened valuation ring to the language. This is made possible via a mixed characteristic resplendency result, as in 
\cite[Section~4]{CHR}.
All proofs in this subsection work exactly as in \cite{CHR}; one just needs to replace
$\RV_{I}$ by $\bigcup_n\RV_{nI}$ everywhere (where $I$ is an ideal of $\cO_K$). For completeness, we nevertheless give most of the details.

Fix a \(\emptyset\)-definable ideal \(I\) of $\cO_K$ which is neither $\{0\}$ nor equal to $\cO_K$.
Call a language $\cL'$ an $\bigcup_n\RV_{nI}$-expansion of $\cL$ if $\cL'$ is obtained from $\cL$ by adding (any) predicates which live on Cartesian products of the (imaginary) definable sets $\RV_{nI}$ for some $n>0$.
(Recall that $\RV_I$ has been defined in \ref{ssec:notation}.)

In the following, we assume that \(K\) is $\kappa$-saturated for some $\kappa > |\cL|$; we call a set ``small'' if its cardinality is less than $\kappa$, and ``large'' otherwise.

\begin{defn}[Mixed-$I$-preparation]
We say that $K$ has \emph{mixed-$I$-preparation} if for every integer \(n\), every set $A\subset K$ and every $(A\cup \RV_{nI})$-definable subset $X\subset K$, there exists a finite $A$-definable set $C \subset K$ and an integer \(m\) such that $X$ is $mI$-prepared by $C$.

We say that $K$ has \emph{resplendent mixed-$I$-preparation} if for every set $A\subset K$, for every $\bigcup_n\RV_{nI}$-expansion $\cL'$ of $\cL$, and
for every $\cL'(A)$-definable subset $X\subset K$,
there exists an integer \(m\) and a finite $\cL(A)$-definable set $C \subset K$ such that $X$ is $mI$-prepared by $C$.
\end{defn}

Clearly, resplendent mixed-$I$-preparation implies mixed-$I$-preparation.

We consider the language $\cL$ as having the sorts \(K\) and \(\RV_{nI}\) for each $n > 0$. Let moreover \(\LI\) be the language with the same sorts, but consisting only of the additive group on \(K\) and the maps \(\rv_{nI}\).

\begin{lem}[Preparation and partial isomorphisms]\label{lem:prep and part iso}
Let $A \leq K$ be a small $\QQ$-sub-vector space. The following are equivalent:
\begin{enumerate}[(i)]
\item Any $\cL(A\cup \bigcup_n \RV_{nI})$-definable set $X\subset K$ can be $mI$-prepared, for some integer \(m\), by some finite set $C\subset A$.
\item For every small subset $A_2\subset K$, $c_1, c_2 \in K$ and all (potentially large) sets $B_1, B_2\subset \bigcup_n\RV_{nI}$ with $\bigcup_n\rv_{nI}(\langle A,c_1\rangle) \subset B_1$,
if $f\colon A B_1 c_1 \to A_2 B_2 c_2$ is a partial $\LI$-isomorphism sending $c_1$ to $c_2$
whose restriction $\restr{f}{AB_1}$ is a partial elementary $\cL$-isomorphism, then the entire $f$
is a partial elementary $\cL$-isomorphism.
\item For all $c_1, c_2 \in K$ and all (potentially large) sets $B\subset \RV_{nI}$ which contain $\bigcup_n\rv_{nI}(\langle A,c_1\rangle)$, any partial $\LI(A\cup B)$-isomorphism $f\colon\{c_1\} \to \{c_2\}$, is a partial elementary $\cL(A\cup B)$-isomorphism.
\end{enumerate}
\end{lem}

\begin{proof}
(i) $\Rightarrow$ (iii):
Let $f$ be as in (iii). We have to check that for every $\cL(A\cup B)$-definable set $X\subset K$, $c_1 \in X$ if and only if $c_2\in X$. By (i), there exists a finite $C \subset A$ and an integer $m$ such that $X$ is $mI$-prepared by $C$. Since $f$ is an $\LI(A\cup B)$-isomorphism and $B$ contains $\rv_{mI}(\langle A,c_1\rangle)$,
for all $a\in C$ and all $r\geq 1$, we have
\[
\rv_{mI}(c_2-a) = \rv_{mI}(f(c_1)-a) = f(\rv_{mI}(c_1-a)) = \rv_{mI}(c_1-a).
\]
Since $X$ is $I$-prepared by $C$, it follows that $c_1 \in X$ if and only if $c_2\in X$.

\medskip

(iii) $\Rightarrow$ (ii):
Let $f$ be as in (ii). Since $f$ is $\cL$-elementary if and only its restriction to every finite domain is, we may assume $B_i$ small. Using the assumption that $\restr{f}{AB_1}$ is $\cL$-elementary, we can extend $(\restr{f}{AB_1})^{-1}$ $\cL$-elementarily to some $g$ defined at $c_2$. Let $c'_1 := g(c_2)$. Then $g\circ f\colon \{c_1\} \to \{c'_1\}$ is a partial $\LI(A\cup B_1)$-isomorphism. Since $\bigcup_n\rv_{nI}(\langle A, c_1\rangle)\subset B_1$, it follows by (iii) that $g\circ f$ is an elementary $\cL$-isomorphism. As $g$ is also $\cL$-elementary, so is $f$.

\medskip

(ii) $\Rightarrow$ (i):
Let $X$ be as in (i), and let $B \subset \bigcup_n\RV_{nI}$ be a finite subset such that $X$ is $\cL(A \cup B)$-definable. Let $\LI'$ denote the expansion of \(\LI\) by the full \(\cL(A)\)-induced structure on \(\bigcup_n\RV_{nI}\).

Consider any $c_1, c_2 \in K$ which have the same qf-$\LI'$-type over $A \cup B$. Then the map $f\colon c_1 \to c_2$ is an $\LI'(A \cup B)$-isomorphism and extends to $f\colon A B_1 c_1 \to A_2 B_2 c_2$, where
$B_i := B \cup \bigcup_n \rv_{nI}(\langle A,c_i\rangle)$. By definition of $\LI'$, the restriction $f|_{AB_1}$ is $\cL$-elementary, so by (ii), the map $f$ is $\cL$-elementary.
Since moreover $f$ is the identity on $A \cup B$, this implies that $c_1$ and $c_2$ have the same $\cL$-type over $A \cup B$.

We just proved that the $\cL(A \cup B)$-type of any element $c \in K$ is implied by its qf-$\LI'(A \cup B)$-type.
By a classical compactness argument (cf. the proof of \cite[Theorem 3.2.5]{TentZiegler}), it follows that any $\cL(A\cup B)$-formula in one valued field variable is equivalent to a quantifier free $\LI'(A\cup B)$-formula. In particular, this applies to our set $X$. Since $\LI'(B)$ is an $\bigcup_n\RV_{nI}$-expansion of $\LI$, \(X\) is indeed \(mI\)-prepared by some finite \(C\subset A\).
\end{proof}

We can now prove, exactly as in 
\cite[Lemma~4.1.16]{CHR}:

\begin{lem}[Back and forth over $\bigcup_n \RV_{nI}$]\label{lem:b-a-f over RV}
Let $K$ have mixed-$I$-preparation. Then the set of partial elementary $\cL$-isomorphisms $f\colon \bigcup_n\RV_{nI}\cup A_1 \to \bigcup_{n}\RV_{nI}\cup A_2$ (where $A_1, A_2$ run over all small subsets of $K$) has the back-and-forth, i.e., given such an $f$ and a $c_1 \in K \setminus A_1$, $f$ can be extended to a partial elementary $\cL$-isomorphism on $\bigcup_n\RV_{nI}\cup A_1 \cup \{c_1\}$.
\end{lem}

\begin{proof}
Let $A_1, A_2, f, c_1$ be given.
We may assume that $\acl_K A_i = A_i$.
Set $B_1 := \bigcup_n\rv_{nI}(\langle A_1, c_1\rangle_\QQ)$. We claim that we can extend $f$ by setting $f(c_1) := c_2$ for any $c_2\in K$ realizing $f_*\qftp_{\LI}(c_1/A_1B_1)$.
To see this, we use \ref{lem:prep and part iso} (i)$\Rightarrow$(ii), i.e., we need to verify that $f|_{A_1B_1}$ is a partial elementary $\cL$-automorphism (this is clear by assumption) and that $f$ is a partial $\LI$-isomorphism. This follows from the restrictions $f|_{A_1B_1c_1}$ and $f|_{\bigcup_n\RV_{nI}}$ being partial $\LI$-isomorphisms, since the only $\LI$-interaction between $A_1c_1$ and $\bigcup_n\RV_{nI}$ is via $B_1$.
\end{proof}

And finally, as in 
\cite[Proposition~4.1.7]{CHR}:

\begin{prop}\label{prop:resp}
Assuming \(K\) is $\kappa$-saturated for some $\kappa > |\cL|$, the following are equivalent:
\begin{enumerate}[(i)]
\item $K$ has mixed-$I$-preparation.
\item $K$ has resplendent mixed-$I$-preparation.
\end{enumerate}
\end{prop}

\begin{proof}
(ii) $\Rightarrow$ (i) is trivial, so we assume (i) and prove (ii). Let $\cL' \supset \cL$ be an $\bigcup_n\RV_{nI}$-expansion and suppose that $X \subset K$ is $\cL'(A)$-definable, for some $A \subset K$. We need to find a finite $\cL(A)$-definable $C \subset K$ $mI$-preparing $X$ for some $m \ge 1$.
We may replace $K$ by a sufficiently saturated elementary extension. We may also assume that $A = \acl_{\cL,K}(A)$ (since every finite $\acl_{\cL,K}(A)$-definable set is contained in a finite $A$-definable set). Now we are in the setting of Lemma~\ref{lem:prep and part iso}.
Applying the lemma for $\cL'$ shows that it suffices to verify that every partial $\LI(\tilde A)$-isomorphism $f\colon \{c_1\} \to \{c_2\}$ is an elementary $\cL'(\tilde A)$-isomorphism, for $\tilde A :=A \cup \bigcup_n\RV_{nI}$, and applying the lemma for $\cL$ shows that such an $f$ is an elementary $\cL(\tilde A)$-isomorphism.
Thus we can finish the proof by showing (more generally) that any elementary $\cL(\tilde A)$-isomorphism $f\colon B_1 \to B_2$ is already an elementary $\cL'(\tilde A)$-isomorphism, for small sets $B_1, B_2 \subset K$. We show that such $f$ preserve every $\cL'(\tilde A)$-formula, by induction over the structure of the formula. For quantifier free formulas, this is clear. To deal with an existence quantifier, use Lemma~\ref{lem:b-a-f over RV} to extend $f$ to a realization of the existential quantifier.
\end{proof}

We now get the analogue of 
\cite[Theorem~4.1.19]{CHR} (about \(\bigcup_n\RV_n\)-expansions preserving Hensel minimality), where the case with $\ell=1$ is replaced by conditions (T1,T2). Once we will have proved \Cref{thm:equivalences:1hmin} in full, we recover the case $\ell=1$ as well.

\begin{prop}\label{prop: RVenrichment mixed}
\begin{itemize}
\item[(i)]
If \(\cT\) is \(\ell\)-\hmix-minimal for some \(\ell \in \{0,\omega\}\), then any \(\bigcup_n\RV_n\)-expansion also is \(\ell\)-\hmix-minimal.

\item[(ii)] If $\cT$ satisfies (T1,T2), then any $(\bigcup_n\RV_{n})$-expansion of $\cT$ also satisfies (T1,T2).
\end{itemize}
\end{prop}

\begin{proof}
(i) follows easily from \Cref{prop:resp} (applied in sufficiently saturated models $K \models \cT$), using $I = \cM_K$ in the case $\ell = 0$, and using $I = B_{<\lambda}(0)$ for every $\lambda \le 1$ in the case $\ell = \omega$. In the latter case, we need to first add $\lambda$ to the language and then get rid of it again using that $\RV$-unions preserve finiteness (\Cref{cor: 1hmixmin RVfiniteness}~(1)), i.e., the full argument is as follows: To $\lambda$-prepare a set $X \subset K$ which is $A \cup \RV_{\lambda}$-definable (where $A \subset K$) in the \(\bigcup_n\RV_n\)-expansion of $K$, apply \Cref{prop:resp} with $I = B_{<\lambda}(0)$ in the language $\cL(\lambda)$ to get a finite $\cL(A,\lambda)$-definable set $C' = \phi(K,\lambda)$ $\lambda$-preparing $X$. By \Cref{cor: 1hmixmin RVfiniteness}~(1), the $\cL(A)$-definable set $C = \bigcup_{\lambda'} \phi(K,\lambda')$ is still finite (and it $\lambda$-prepares $X$ since it contains $C'$).

For (ii), let $\cL'$ be an $(\bigcup_n\RV_{|n|})$-expansion of the language $\cL$, and let $K$ be a model of the expansion of $\cT$.
First of all, by part (i) (and Lemma~\ref{lem: property 3 basic props} (2)) we know that $\Th_{\cL'}(K)$ is still $0^\mix$-h-minimal. Let $f: K\to K$ be $\cL'(A)$-definable, for some $A\subset K\cup \RV_n$. By Proposition~\ref{prop:resp}, the $\cL'(A \cup \{x\})$-definable set $\{f(x)\}$ can be prepared by a finite $\cL(A \cup \{x\})$-definable set $C_x$; note that this implies $f(x) \in C_x$. We may assume (using compactness) that
$C_x$ is definable uniformly in $x$. By Lemma \ref{lem: basics 1hmixmin}, we find an $\cL(A)$-definable family of injective maps $g_x: C_x\to \RV_m^k$. We define $h: K\to \RV_m^k: x\mapsto g_x(f(x))$ and
\[
\tilde{f}: K\times \RV_m^k\to K: (x, \zeta)\mapsto \begin{cases} g_x^{-1}(\zeta) & \text{ if } \zeta\in g_x(C_x), \\
	0 & \text{ else.} \end{cases}
\]
Note that $f(x) = \tilde{f}(x, h(x))$ for every $x\in K$, and that $h$ is $\cL'(A)$-definable while $\tilde{f}$ is $\cL(A)$-definable. By Lemma \ref{lem: basics 1hmixmin}, we can find an $\cL(A)$-definable family of finite sets $C_\xi$, for $\xi\in \RV_m^k$ such that the map $x\mapsto \tilde{f}(x, \xi)$ satisfies (T1$^\mix$) with integer $p$ and (T2) with respect to $C_\xi$. Let $C_1$ be the union of all of these $C_\xi$, which is still finite by Lemma \ref{lem: property 3 basic props}. By $0^\mix$-h-minimality, there exists a finite $\cL'(A)$-definable set $C_2\subset K$ and a positive integer $p'$ such that $h$ is constant on the balls $|p'|$-next to $C_2$. Now take $C = C_1\cup C_2$ and let $q = pp'$. Then $C$ is a finite $\cL'(A)$-definable set for which $f$ satisfies properties (T1$^\mix$) and (T2) with respect to $C$ (with integer $q$).
\end{proof}

\begin{remark}\label{rem: RVenrichment mixed}
After proving \Cref{thm:equivalences:1hmin} (which will be done in the next subsection), the resplendency statement (i) of \Cref{prop: RVenrichment mixed} follows for $\ell=1$. Indeed, this is direct from (ii) of \Cref{prop: RVenrichment mixed} and \Cref{thm:equivalences:1hmin}. Preservation of $\ell$-h-minimality under $\RV$-enrichment is generalized in \cite[Theorem 2.2.5]{Verm:h-min} to all $\ell\geq 0$ in equi-characteristic zero.
\end{remark}

\subsection{Proof of Theorem \ref{thm:equivalences:1hmin}}\label{sec:equiv:proof}

We are now ready to complete the proof of Theorem \ref{thm:equivalences:1hmin}.

\begin{proof}[Proof of Theorem \ref{thm:equivalences:1hmin}]
\textbf{(1) implies (2).} This is \Cref{cor:prep:mix}.

\textbf{(2) implies (3).} This is \Cref{cor.T1.T2}.

\textbf{(3) implies (4).} Let $|\cdot |_c$ be a nontrivial equicharacteristic zero coarsening of the valuation $|\cdot|$.

The coarsened valuation ring $\cO_{K,c}$ is a pullback of some subset of $\RV$. Hence by Item (ii) of \Cref{prop: RVenrichment mixed}, the theory of $K$ in $\cL' = \cL\cup \{\cO_{K,c}\}$ still satisfies (T1,T2). We work in this language $\cL'$.

We will use 
\cite[Theorem~2.9.1]{CHR}, which says that the theory of $K$ in $\cL_c$ with respect to the norm $|\cdot |_c$ is $1$-h-minimal if for every set $A\subset K\cup \RV_c$ and every $A$-definable function $f: K\to K$ the following hold:
\begin{enumerate}
\item there exists some finite $A$-definable set $C$ such that for every ball $B$ which is $|1|_c$-next to $C$ there is a $\mu_B\in \Gamma_{K,c}$ such that for $x,y\in B$
\[
|f(x)-f(y)|_c = \mu_B |x-y|_c, \text{    and}
\]
\item the set $\{y\in K\mid f^{-1}(y) \text{ is infinite}\}$ is finite.
\end{enumerate}
Let us say that a finite set $C$ \emph{$\mathfrak{m}_c$-prepares a function $f$} if (1) holds for $f$ and $C$.

So let $f: K\to K$ be $A$-definable for some $A\subset   K\cup \RV_{c}$. By adding parameters from $K$ to the language we may as well assume that $A\subset   \RV_c$. We will find a finite $\emptyset$-definable $C$ which $\mathfrak{m}_c$-prepares $f$. We induct on $r=\# A$. The case $r = 0$ is clear. Indeed, in that case there even exists a finite $\cL$-definable set $C$ such that $f$ satisfies (T1) with respect to $C,m$. In particular $C$ will $\mathfrak{m}_c$-prepare $f$.

So suppose that $r>0$. Let $A=A'\cup \{\zeta\}$ and consider an $A'$-definable family $(f_y)_y$ for $y\in K$ with $f_y = f$ if $\rv_c y=\zeta$. Here we use the fact that $\mathfrak{m}_c$ is $\emptyset$-definable. Let $Y = \rv_c^{-1}(\zeta)$. By induction on $r$ there exist for every $y\in K$ a finite $y$-definable $C_y$ $\mathfrak{m}_{c}$-preparing $f_y$. Fix a $y$-definable injection $h_y: C_y\to \RV_{|n_0|}^k$ for some integer $n_0$. Such a map exists by Lemma \ref{lem: basics 1hmixmin}(4) which holds under (T1,T2) by Lemma \ref{lem: property 3 basic props}.
By compactness we may assume that both $C_y$ and $h_y$ form $\emptyset$-definable families with a parameter $y$. For $\eta\in \RV_{n_0}^k$ consider
\begin{align*}
Y_\eta &= \{y\in K\mid \eta\in \operatorname{im}(h_y)\} \\
g_\eta&: Y_\eta\to K: y\mapsto h_{y}^{-1}(\eta).
\end{align*}
Both of these are $\eta$-definable and via compactness we may assume that the entire family is $\emptyset$-definable with a parameter. Note that
\[
\bigcup_{y\in K} \{y\}\times C_y = \bigsqcup_{\eta\in \RV_{n_0}^k} \operatorname{graph} (g_\eta).
\]
Let $D_\eta$ be a finite $\eta$-definable set $|m_0|$-preparing $Y_\eta$ and such that $g_\eta$ (extended by zero) satisfies Lemma \ref{lem:property (3) balls map to balls} with respect to $D_\eta$ and $m_0$. This is possible by Lemma \ref{lem:property (3) balls map to balls}. Moreover take $D_\eta$ to $\mathfrak{m}_c$-prepare the $(A'\cup\{\eta\})$-definable
\[
\tilde{f}_\eta: K\to K: y\mapsto f_y(g_\eta(y)).
\]
This is via induction on $r$ together with the fact that we can add $\eta$ to the language and preserve (T1,T2). By compactness, we have without loss that $D_\eta$ forms a $\emptyset$-definable family, with a parameter $\eta$, and that $m_0$ does not depend on $\eta$. Now put $D=\{0\}\cup \bigcup_{\eta\in \RV_{|n_0|}^k} D_\eta$, which is a finite $\emptyset$-definable set by Lemma \ref{lem: basics 1hmixmin}.

If $D\cap Y\neq \emptyset$ then $C = \bigcup_{y\in D} C_y$ is a finite $\emptyset$-definable set which $\mathfrak{m}_c$-prepares $f$. So suppose that $D\cap Y = \emptyset$. Since $0\in D$ and as $Y$ is an $\rv_c$-fibre, we have that $Y$ is $\mathfrak{m}_c$-next to $D$.

Use Lemma \ref{lem: basics 1hmixmin} to find integers $p,k'$ and a $\emptyset$-definable map
\[
\chi: K\to \RV_{|p|}^{k'}
\]
such that any ball $|m_0|$-next to $D$ is a union of fibres of $\chi$. Use \Cref{cor: 1hmixmin preparing families} to find a finite $\emptyset$-definable set $C$ and a positive integer $n'$ such that $C$ $|n'|$-prepares every set $g_\eta(\chi^{-1}(\xi))$ for any $\eta\in \RV_{|n_0|}^k$ and any $\xi\in \RV_{|p|}^{k'}$ for which $\chi^{-1}(\xi)\subset Y_\eta$. Then we have that for any ball $B$ $|m_0|$-next to $D$ with $B\subset Y_\eta$ that $g_\eta(B)$ is $|n'|$-prepared by $C$. We will prove that $C$ $\mathfrak{m}_c$-prepares $f$.

We claim that $C$ $\mathfrak{m}_c$-prepares $g_\eta(Y)$ for any $\eta\in \RV_{n_0}^k$ with $Y\subset Y_\eta$. To this end, we will show that for any $\lambda\in \Gamma_K^\times, \lambda\leq 1$, if $B_0$ is a ball $|m_0|\lambda$-next to $D$ then $g_\eta(B_0)$ is $|m_0^2n'|\lambda$-prepared by $C$ whenever $B\subset Y_\eta$. So let $B_0$ be $|m_0|\lambda$-next to $D$ for some $\lambda\in \Gamma_K^\times, \lambda \leq 1$ and let $B$ be the ball $|m_0|$-next to $D$ containing $B_0$. Suppose that $B\subset Y_\eta$ and let $\mu_B$ be as in (T1,T2) for $g_\eta$ on $B$. By Lemma \ref{lem:property (3) balls map to balls}, for any $x\in B_0$ the open ball of radius $|m_0|\mu_B\radop(B_0)$ around $g_\eta(x)$ is contained in $g_\eta(B)$. By the same Lemma, $g_\eta(B)$ is contained in an open ball of radius $|1/m_0|\mu_B\radop(B)$. Now let $B'_0$ be the ball $|m_0^2n'|\lambda$-next to $C$ containing $g_\eta(x)$ and let $B'$ be the ball $|n'|$-next to $C$ containing $B'_0$. Since $g_\eta(B)$ is $|n'|$-prepared by $C$, $B'\subset g_\eta(B)$. Hence
\begin{align*}
\radop B_0' = |m_0^2|\lambda\radop B' \leq |m_0|\lambda \mu_B \radop B = |m_0|\mu_B\radop B_0.
\end{align*}
Therefore $B'_0\subset g_\eta(B_0)$ and we conclude that $C$ $|m_0^2n'|\lambda$-prepares $g_\eta(B)$.

We now prove the claim that $C$ $\mathfrak{m}_c$-prepares $g_\eta(Y)$ for any $\eta\in \RV_{n_0}^k$ with $Y\subset Y_\eta$. Fix such an $\eta$ and recall that $Y$ is $\mathfrak{m}_c$-next to $D$. Let $B_0$ be the ball $|m_0|$-next to $D$ containing $Y$. By property T1, $g_\eta$ is either injective or constant on $B_0$. If $g_\eta$ is constant on $B_0$ there is nothing to show, so assume that $g_\eta$ is injective on $B_0$. Let $B'$ be a ball $\mathfrak{m}_c$-next to $C$. If $B'\cap g_\eta(Y)=\emptyset$ then there is nothing to show. So we may assume that $g_\eta(y)$ and $B'$ are not disjoint. Assume that the claim is false, so there exists some $z\in B'\setminus g_\eta(Y)$. Now $B'$ is contained in a ball $|n|$-next to $C$, which is in turn contained in $g_\eta(B_0)$ by construction of $C$. By injectivity of $g_\eta$ on $B_0$, there is a unique $x\in B_0$ with $g_\eta(x)=z$. We necessarily have that $x\notin Y$. Take an open ball $B_1$ around $Y$ not containing $x$. This ball is $|m_0|\lambda$-next to $D$ for some $\lambda\leq 1$ in $\Gamma_K^\times$. Since $Y$ is $\mathfrak{m}_c$-next to $D$ we have that the image of $|m_0|\lambda$ is $1$ under the map $\Gamma_K^\times \to \Gamma_{K,c}^\times$. By the previous paragraph, $g_\eta(B_1)$ is $|m_0^2n'|\lambda$-prepared by $C$. Now similarly, $|m_0^2n'|\lambda$ becomes $1$ in $\Gamma_{K,c}^\times$ and so $B'\subset g_\eta(B_1)$. But then $z\in g_\eta(B_1)$, contrary to assumption.

We now prove that $C$ $\mathfrak{m}_c$-prepares $f$. Suppose the result is false, and let $B$ be a ball $\mathfrak{m}_c$-next to $C$ on which the result does not hold. So there exist distinct $x,x',x''\in B$ and $\mu_1, \mu_2\in \Gamma_{K,c}$ with $\mu_1\neq \mu_2$ and
\begin{align*}
|f(x)-f(x')|_c &= \mu_1|x-x'|_c \\
|f(x)-f(x'')|_c &= \mu_2|x-x''|_c.
\end{align*}
Assume that $|x-x'|_c\geq |x-x''|_c$ and let $B_1$ be the smallest closed ball containing $x,x'$ for $|\cdot |_c$. Fix a $y_0\in Y$. Since $C_{y_0}$ $\mathfrak{m}_c$-prepares $f_{y_0}=f$ and since $x,x',x''\in B_1$ we have $C_{y_0}\cap B_1\neq \emptyset$. So let $\eta\in \RV^k$ be such that $g_\eta(y_0)\in B_1$. Recall that $Y$ is a maximal ball disjoint from $D$ for the coarsened norm. By our construction of $D$, $g_\eta$ is defined on all of $Y$. We claim that since $D$ $\mathfrak{m}_c$-prepares $\tilde{f}_\eta$, $g_\eta(Y)$ cannot contain all of $B_1$. Indeed, suppose to the contrary that $g_\eta(y)=x, g_\eta(y')=x'$ and $g_\eta(y'')=x''$ for some $y,y',y''\in Y$. Then there are $\mu_{g_\eta}, \mu_{\tilde{f}_\eta}\in \Gamma_{K,c}$ such that
\begin{align*}
|f(x)-f(x')|_c &= \mu_1|x-x'|_c = \mu_1\mu_{g_\eta}|y-y'|_c \\
&= |\tilde{f}_\eta(y)-\tilde{f}_\eta(y')|_c = \mu_{\tilde{f}_\eta}|y-y'|_c.
\end{align*}
Hence $\mu_1 = \mu_{\tilde{f}_\eta}/\mu_{g_\eta}$. But exactly the same reasoning works for $\mu_2$, contrary to our assumption that $\mu_1$ and $\mu_2$ are distinct. Thus $g_\eta(Y)$ cannot contain all of $B_1$.

Therefore
\[
\radopeqc(g_\eta(Y))\leq \radcleqc(B_1) = |x-x'|_c.
\]
Now $C$ $\mathfrak{m}_c$-prepares $g_\eta(Y)$ so there exists some $c\in C$ with
\[
|g_\eta(y_0)-c|_c \leq \radopeqc(g_\eta(Y))\leq |x-x'|_c.
\]
So we have
\[
|x-c|_c \leq \max\{|x-g_\eta(y_0)|_c, |g_\eta(y_0) - c|_c \} \leq |x-x'|_c.
\]
But this would mean that $x$ and $x'$ are not in the same ball $\mathfrak{m}_c$-next to $C$, contradiction.

\textbf{(4) implies (1).} Trivial.
\end{proof}

\section{Geometry of definable sets in 1-h-minimal  structures}\label{ssec:coarsenings}

In this section
we develop our main geometrical results under 1-h-minimality, often using the equivalences of \Cref{thm:equivalences:1hmin}.

\subsection{Taylor approximation}\label{sec:Tay}

Results which simply state that every definable set or function of a certain kind has some nice (language-independent) properties immediately follow for valued fields of mixed characteristic via coarsenings; in particular, we have a good dimension theory, similar to 
\cite[Section~5.3]{CHR}. Here the dimension of a (non-empty) definable set $X\subset K^n$ is the maximal integer $m$ such that there is a $K$-linear function $\ell: K^n\to K^m$ for which $\ell(X)$ has non-empty interior. We denote it by $\dim X$ and put $\dim \emptyset = -\infty$.

 More generally we have the following:
\begin{prop}[Language-independent properties]\label{prop:eqc-easy}
Let $\cT$ be a $1$-h-minimal theory. Then the conclusions of the following results hold for any model $K$ of $\cT$:

\begin{enumerate}
\item (Almost everywhere continuity 
\cite[Theorem~5.1.1]{CHR}). For every definable function $f: X\subset K^n \to K$, the set of $u\in X$ at which $f$ is continuous in a neighbourhood of $u$ is dense in $X$.
\item (Almost everywhere $C^k$ 
\cite[Theorem~5.1.5]{CHR}). For every definable function $f: X\subset K^n \to K$, the set of $u\in X$ at which $f$ is $C^k$ in a neighbourhood of $u$ is dense in $X$.
\item (Dimension theory 
\cite[Proposition~5.3.4]{CHR}). Let $X,Y\subset K^n, Z\subset K^m$ be non-empty and definable, and let $f: X\to Z$ be definable. Then the following hold
\begin{enumerate}
\item We have $\dim X\geq d$ if and only if there is a coordinate projection $\pi: K^n\to K^d$ such that $\pi(X)$ has non-empty interior.
\item $\dim (X\cup Y) = \max \{\dim X, \dim Y\}$.
\item For any integer $d$, the set of $z\in Z$ such that $\dim f^{-1}(z) = d$ is definable over the same parameters of $f$.
\item If all fibres of $f$ have dimension $d$, then $\dim X = d + \dim Z$.
\item There exists an $x\in X$ such that for every open ball $B\subset K^n$ around $x$, $\dim X = \dim (X\cap B)$.
\item We have that $\dim (\overline{X} \setminus X) < \dim X$ where $\overline{X}$ is the closure of $X$.
\end{enumerate}
\end{enumerate}
\end{prop}

\begin{proof}
We may work in a model $K \models \cT$ such that $|\cdot|_{\eqc}$ is non-trivial. Then every $\cL$-definable object is in particular $\cL_{\eqc}$-definable, so all of the above equi-characteristic $0$ results apply to the definable objects in question and yield the desired mixed-characteristic result, except in the case of Item~3.c.
(Concerning 
\cite[Theorem~5.1.1]{CHR} and 
\cite[Theorem~5.1.5]{CHR}, note that $|\cdot|$ and $|\cdot|_{\eqc}$
induce the same topology and hence equivalent notions of continuity and derivatives.)

Now Item~3.c (definability of dimension) can easily be reproved directly in $\cL$, using \Cref{lem: basics 1hmixmin} and that $0$-dimensional is equivalent to finite.
\end{proof}

Recall that if $f: K\to K$ is a function which is $C^r$ at a point $x_0\in K$, then we define the order $r$ Taylor polynomial of $f$ at $x_0$ as
\[
T_{f,x_0}^{\leq r}(x) = \sum_{i=0}^{r} \frac{f^{(i)}(x_0)}{i!}(x-x_0)^i.
\]
Now we find back the precise analogue of 
\cite[Theorem~3.2.2]{CHR} on Taylor approximation.

\begin{thm}[Taylor approximation]\label{cor:high-ord:mix}
Suppose that  $\cT$ is $1$-h-minimal and let $K$ be a model of $\cT$. Let $f\colon K \to K$ be an $A$-definable function, for $A\subset K\cup \RV_{|n|}$, and let $r \in \NN$ be given. Then there exists a finite $A\cap K$-definable set $C$ and an integer $m \ge 1$ such that for every ball $B$ $|m|$-next to $C$, $f$ is $(r+1)$-fold differentiable on $B$,
$|f^{(r+1)}|$ is constant on $B$,
and we have:
\begin{equation}\label{eq:t-higher:mix}
 |f(x) -  T^{\le r}_{f,x_0}(x) | = \left|\frac{1}{(r+1)!} f^{(r+1)}(x_0)\cdot (x-x_0)^{r+1}\right|
\end{equation}
for every $x_0, x \in B$.
\end{thm}

\begin{proof}
As usual, we may simply assume that $A=\emptyset$, since we can add parameters from $K\cup\RV_{|n|}$ to $\cL$ and preserve 1-h-minimality. We first prove a slightly weaker version of this result, where the conclusion states that
\begin{equation}\label{eq:weak-t-higher:mix}
|f(x) -  T^{\le r}_{f,x_0}(x) | \leq \left|\frac{1}{m}\cdot f^{(r+1)}(x_0)\cdot (x-x_0)^{r+1}\right|,
\end{equation}
for some positive integer $m$.
Proceeding as in the proof of Corollary~\ref{cor:prep:mix}, we assume that $K$ is sufficiently saturated and sufficiently homogeneous and we use 
\cite[Theorem~3.2.2]{CHR} and Lemma~\ref{lem:mix} to find a finite $\cL$-definable set $C$ such that $f$ is differentiable away from $C$, and such that for $x_0, x$ in the same ball $|1|_\eqc$-next to $C$, we have
\begin{equation}\label{eq:t-higher:coarse}
|f(x) -  T^{\le r}_{f,x_0}(x) |_\eqc \leq |f^{(r+1)}(x_0)\cdot (x-x_0)^{r+1}|_\eqc.
\end{equation}
Now suppose that there exists no $m \ge 1$ satisfying the condition involving (\ref{eq:weak-t-higher:mix}). As before, we then use $\aleph_0$-saturation to find a pair $(x_0, x) \in K^2$ of points which lie in the same ball $|1|_\eqc$-next to $C$ and such that (\ref{eq:t-higher:mix}) fails for every $m$. The latter means that (\ref{eq:t-higher:coarse}) fails for $x_0, x$, so we have a contradiction to our choice of $C$.

For the weakened theorem, it remains to ensure that $|f^{(r+1)}|$ is constant on balls $|m|$-next to $C$. By applying \Cref{cor: 1hmixmin preparing families} to the graph of $x \mapsto \rv(f^{(r+1)}(x))$, we can enlarge $C$ so that this works for balls $|m'|$-next to $C$, for some $m' \ge 1$. Now the weaker statement holds using $m \cdot m'$ (since $|m \cdot m'| \le |m|, |m'|$).

We now explain how to get the statement with $(r+1)!$ instead of $m$. We can find a finite $A\cap K$-definable $C$ and an integer $m$ such that $f$ satisfies the weakened version of this theorem with order $r+1$. We now enlarge $C$ and find an integer $n$ such that $f^{(r+1)}$ satisfies the conclusion of Lemma~\ref{lem:property (3) balls map to balls} with respect to $C$ and $n$. We may moreover assume by Lemma~\ref{prop: 1hmixmin preparing families} that $\rv\circ f^{(r+1)}$ is constant on balls $|n|$-next to $C$. Then there exists an integer $n'$ such that for any ball $B$ $|n|$-next to $C$ and any $x_0\in B$ we have that
\[
|f^{(r+2)}(x_0)|\radop B \leq \frac{1}{|n'|}|f^{(r+1)}(x_0)|.
\]
Now at the cost of enlarging $n$ (and making $B$ smaller) we obtain
\[
|f^{(r+2)}(x_0)|\radop B \leq |m||f^{(r+1)}(x_0)|.
\]
Thus we obtain for $x, x_0\in B$ that
\[
|f(x)-T_{f, x_0}^{\leq r+1}(x)|\leq \left|\frac{1}{m}f^{(r+2)}(x_0)\cdot (x-x_0)^{r+2}\right| \leq \frac{|f^{(r+1)}(x_0)(x-x_0)^{r+2}|}{\radop(B)}.
\]
But since $|x-x_0|< \radop(B)$ we get that
\[
|f(x)-T_{f, x_0}^{\leq r+1}(x)|< |f^{(r+1)}(x_0)(x-x_0)^{r+1}|.
\]
Now apply the ultrametric triangle inequality to get
\begin{align*}
|f(x)-T_{f, x_0}^{\leq r}(x)| & \le \max\{ |f(x) - T_{f, x_0}^{\leq r+1}(x)|, \left| \frac{f^{(r+1)}(x_0)}{(r+1)!}(x-x_0)^{r+1}\right|\} \\
&= \left| \frac{1}{(r+1)!}\right| \cdot |f^{(r+1)}(x_0)(x-x_0)^{r+1}|.\qedhere
\end{align*}

\end{proof}

Surprisingly, this implies an exact version of property (T1$^\mix$) and of Lemma \ref{lem: 1hmix valjac}, namely, an exact valuative Jacobian property similar to 
\cite[Corollary~3.1.6]{CHR}, and, we get an ever finer version with $\rv$ in \Cref{cor.jacprop}.

\begin{cor}[Valuative Jacobian property]\label{cor:T1exactmix}
Let $\cT$ be a 1-h-minimal theory. Let $K$ be a model of $\cT$ and let $f: K\to K$ be $A$-definable, for some $A\subset   K\cup \RV_{|n|}$. Then there exists a finite $\emptyset$-definable set $C$ and an integer $m$ such that for any ball $B$ $|m|$-next to $C$, $f$ is differentiable on $B$, $|f'|$ is constant on $B$, and,
for all $x, x'\in B$ one has
\[
|f(x)-f(x')| = |f'(x)|\cdot|x-x'|.
\]
Moreover, we can ensure that if $B'\subset B$ is an open ball then $f(B')$ is either a singleton or an open ball of radius
\[
|f'(x)|\cdot\radop B'.
\]
\end{cor}
\begin{proof}
Apply \Cref{cor:high-ord:mix} for $r=0$ for the first statement. For the moreover part, follow the proof of 
\cite[Lemma~2.9.6]{CHR}.
\end{proof}

The following finer variant with $\rv$ instead of the valuation is useful for motivic integration, for the change of variables formulas and for Fourier transforms.

\begin{cor}[Jacobian property]\label{cor.jacprop}
Let $\cT$ be a 1-h-minimal theory. Let $K$ be a model of $\cT$ and let $f: K\to K$ be $A$-definable, for some $A\subset   K\cup \RV_{|n|}$. Then there exists a finite $A$-definable set $C$ and an integer $m$ such that for every $\lambda\in \Gamma_K^\times$, $\lambda\leq 1$ and every ball $|m|\lambda$-next to $C$ we have:
\begin{enumerate}
\item The derivative $f'$ exists on $B$ and $\rv_\lambda\circ f'$ is constant on $B$.
\item For any $x, x'\in B$ we have $\rv_\lambda \left( \frac{f(x)-f(x')}{x-x'}\right)= \rv_\lambda (f')$.
\item If $B'\subset B$ is any open ball, then $f(B')$ is either a point or an open ball.
\end{enumerate}
\end{cor}
\begin{proof}
Identical to 
\cite[Corollary~3.2.7]{CHR} (whose proof is below \cite[Definition~3.2.8]{CHR}).
\end{proof}

\subsection{Algebraic Skolem functions}

We show that as in the equicharacteristic zero case, one can obtain $\acl=\dcl$ by an $\bigcup_n\RV_n$-expansion. Recall that a theory satisfies $\acl=\dcl$ if for for any model $K$ and any subset $A\subset K$ we have $\acl_K(A) = \dcl_K(A)$. Equivalently, this means that algebraic Skolem functions exist in any model of the theory. The following is the analogue of 
\cite[Proposition~2.6.12]{CHR}.

\begin{lem}\label{lem:stable.embed}
Suppose that $\Th(K)$ is $0$-h-minimal (possibly of mixed characteristic) and let $A\subset K$. Let $X\subset \RV^k_{|m|}$ be an $A$-definable set, for some integers $k,m$. Then $X$ is $\dcl_{\cup_n \RV_{|n|}}(A)$-definable.
\end{lem}
\begin{proof}
We can assume that $A$ is finite, and induct on $\#A$. The case where $A=\emptyset$ is clear, so write $A = \hat{A}\cup \{a\}$. Then we find an $\hat{A}$-definable family
\[
Y\subset K\times \RV_{|m|}^k
\]
such that $Y_a = X$. Using Proposition~\ref{prop: 1hmixmin preparing families} we find a finite $\hat{A}$-definable set $C$ and an integer $p$ such that either $a\in C$, or for every $a'\in K$ in the same ball $|p|$-next to $C$ we have $Y_{a'} = Y_a = X$. In both cases, we use Lemma~\ref{lem: basics 1hmixmin}(5) to find an $\hat{A}$-definable map $f: K\to \RV_{|q|}^{k'}$ for some integers $q,k'$ such that any ball $|p|$-next to $C$ is a union of fibres of $f$. Then $Y_a = X$ is definable over $\hat{A}\cup\{f(a)\}$, so there is a family $Z\subset \RV_{|m|}^{k}\times \RV_{|q|}^{k'}$ such that $X = Z_{f(a)}$. By induction $Z$ is definable over $\dcl_{\cup_n \RV_{|n|}}(\hat{A})$. Since $f(a)\in \RV_{|q|}^{k'}$ is definable over $A$, we find that $X$ is definable over $\dcl_{\cup_n \RV_{|n|}}(A)$.
\end{proof}

\begin{lem}\label{lem:acl}
Suppose that $\Th(K)$ is $0$-h-minimal (possibly of mixed characteristic) and that for any integer $m$ and any set $A\subset \RV_{|m|}$ we have $\acl_{\cup_n\RV_{|n|}}(A) = \dcl_{\cup_n\RV_{|n|}}(A)$. Then for any set $A\subset K$ we have $\acl_K(A) = \dcl_K(A)$.
\end{lem}
\begin{proof}
See also 
\cite[Lemma~4.3.1]{CHR}. Let $C\subset K$ be a finite $A$-definable set. By Lemma~\ref{lem: basics 1hmixmin}(4) there exists an $A$-definable bijection $f: C\to C'\subset \RV_{|n|}^k$ for some integers $n,k$. By the previous lemma, $C'$ is already definable over $\dcl_{\cup_p \RV_{|p|}}(A)$. So by the assumption $C'$ is contained in $\dcl_{\RV_{|n|}}(A)$, and by pulling back via $f$ we have that $C\subset \dcl_K (A)$.

\end{proof}

\begin{prop}\label{prop:acl}
Suppose that $\Th(K)$ is $1$-h-minimal of possibly mixed characteristic in a language $\cL$. Then there exists a $\cup_n \RV_n$-expansion $\cL'$ of $\cL$ such that $\Th_{\cL'}(K)$ is still $1$-h-minimal and for any model $K'$ of $\Th_{\cL'}(K)$ and any subset $A\subset K'$ we have
\[
\acl_{\cL',K'}(A) = \dcl_{\cL',K'}(A).
\]
\end{prop}

\begin{proof}
This is identical to 
\cite[Proposition~4.3.3]{CHR}. We simply let $\cL'$ be the extension of $\cL$ by a predicate for every subset of $\RV_{|n|}^k$ (for all integers $n,k$). The theory remains $1$-h-minimal in view of Proposition~\ref{prop: RVenrichment mixed}. By construction we moreover have for every model $K'$ of $\Th_{\cL'}(K)$ and every subset $A\subset \RV_{|m|}^k$ that $\acl_{\cup_n \RV_{|n|}}(A) = \dcl_{\cup_n \RV_{|n|}}(A)$, so that we can conclude by the previous lemma.
\end{proof}

By an $\cup_n \RV_n$-expansion of the language as in \Cref{prop:acl} we can thus always obtain $\acl=\dcl$. The following lemma provides a tool for returning to the original language.

\begin{lem}[Undoing $\RV$-expansions]\label{lem:undo-K-to-RV:mixed}
Suppose that $\Th_{\cL}(K)$ is 0-h-minimal and that $\cL'$ is an $\cup_n\RV_{|n|}$-expansion of $\cL$. Let $\chi'\colon K^n \to \RV^{k'}_{|m'|}$ be an $\cL'$-definable map (for some positive integers $k' $ and $m'$).
Then there exists an $\cL$-definable map $\chi\colon K^n \to \RV^k_{|m|}$ (for some $k$ and $m$) such that $\chi'$ factors over $\chi$, i.e. $\chi' = g\circ \chi$,
for some function $g\colon \RV_{|m|}^k \to \RV_{|m'|}^{k'}$ (which is automatically $\cL'$-definable).
\end{lem}
\begin{proof}
The proof is similar to 
\cite[Lemma~4.3.4]{CHR}.

We induct on $n$. The case $n=0$ is clear, so let $n>0$ and for $a\in K$ define $\chi'_a: K^{n-1}\to \RV^{k'}_{|m'|}: x\mapsto \chi'(a,x)$. By induction and compactness, we find an $\cL$-definable family of maps $\chi_a: K^{n-1}\to \RV^\ell_{|m|}$ with a parameter $a\in K$ such that for every $a\in K$ we have $\chi'_a = g_a\circ \chi_a$ for some $\cL'$-definable family of maps $g_a$. Now consider the set
\[
W = \{(a, \zeta', g_a(\zeta'))\mid a\in K, \zeta'\in \RV^k_{|m|}\}
\]
and use Lemma~\ref{prop: 1hmixmin preparing families} to uniformly $|p|$-prepare $W$ with some finite $\cL'$-definable set $C'$. Using Proposition~\ref{prop:resp} we find a finite $\cL$-definable set $C$ containing $C'$. By Lemma~\ref{lem: basics 1hmixmin}(5), we find $f: K\to \RV_{|q|}^{\ell'}$ an $\cL$-definable map such that every ball $|p|$-next to $C$ is a union of fibres of $f$. We can assume that $|q|\leq |m|$ and so by lifting $\chi_a$ to a map $K^{n-1}\to \RV_{|q|}^\ell$, the $\cL$-definable map $\chi(a,x) = (f(a), \chi_a(x))$ will be as desired.
\end{proof}

\subsection{Cell decomposition and (higher dimensional) Jacobian property}\label{sec:cd.mix}

As in 
\cite[Section~5.2]{CHR}, we approach the cell decomposition results in a simpler way than usual (compared to e.g.~\cite{Pas}), by imposing $\acl=\dcl$ in $K$, or in other words, that we have algebraic Skolem functions in $K$. This is harmless by \Cref{prop:acl} and \Cref{lem:undo-K-to-RV:mixed}.
In fact, the assumption that $\acl$ equals $\dcl$ does more than just simplify the arguments: it also allows one to formulate stronger results like piecewise Lipschitz continuity results as in Theorem \ref{thm:cd:alg:piece:Lipschitz.mix}.
This is similar to \cite{CFL}, where this condition is also used for similar reasons.

In the following, for $m \le n$, we denote the projection $K^n \to K^{m}$ to the first $m$ coordinates by $\pi_{\le m}$, or also by $\pi_{<m+1}$.

\begin{defn}[Cells, twisted boxes]\label{defn:cell.mix}
Fix any parameter set $A \subset K\eq$ and consider a non-empty $A$-definable set $X\subset  K^n$ for some $n$, an integer $N$, and, for $i=1,\ldots, n$, values $j_i$ in $\{0,1\}$ and $A$-definable functions $c_i:\pi_{<i}(X)\to K$.
Then $X$ is called an \emph{$A$-definable cell} of depth $|N|$, with \emph{center tuple $c=(c_i)_{i=1}^n$} and of \emph{cell-type $j= (j_i)_{i=1}^n$} if it is of the form
$$
X = \{x\in K^n\mid  (\rv_{|N|}(x_i-c_i(x_{<i})))_{i=1}^n
   \in R  \},
$$
for a (automatically $A$-definable) set
$$
R \subset  \prod_{i=1}^n (j_i\cdot \RV_{|N|}^\times),
$$
where $x_{<i}=\pi_{<i}(x)$ and where
$0 \cdot \RV_{|N|}^\times = \{0\} \subset \RV_{|N|}$, and $1 \cdot \RV_{|N|}^\times = \RV_{|N|}\setminus \{0\}$.

If $X$ is such a cell, then, for any $r\in R$, the subset
$$
\{x \in K^n\mid  \rv_{|N|}(x_i-c_i(x_{<i})))_{i=1}^n   = r  \},
$$
of $X$ is called a \emph{twisted box} of the cell $X$.
We also call $X$ itself a twisted box if it is a cell consisting
of a single twisted box (i.e., if $R$ is a singleton).
\end{defn}

Note that a cell of depth $|N|$ is automatically also a cell of depth $|M|$ for any $|M|\leq |N|$, by taking the preimage of the set $R$ under the projection map $\RV_{|M|}^n\to \RV_{|N|}^n$.

\begin{example}
Consider the set of non-zero cubes $X$ in $\QQ_3$ from Example~\ref{eg:cubes.Q3}. This is a $(1)$-cell of depth $|3|$ with cell center $0$. Indeed, $X$ is a union of fibres of the map $\rv_{|3|}: \QQ_3\to \RV_{|3|}$.

For a two-dimensional example, let $Y\subset \QQ_3^2$ be the set of $(x,y)\in \ZZ_3^2$ for which $y-x$ is a non-zero cube. This is a $(1,1)$-cell of depth $|3|$ with cell tuple $c_1 = 0$ and $c_2(x) = x$. In more detail, let $R \subset \RV_{|3|}$ be the set $\rv_{|3|}(X)$ of $\rv_{|3|}$-values of cubes in $\QQ_3$. Then we have that
\[
Y = \{(x,y)\in \ZZ_3^2 \mid (\rv_{|3|}(x - 0), \rv_{|3|}(y - x))\in \RV_{|3|}\times R\}.
\]
\end{example}

As in 
\cite[Section~5.2]{CHR}, we first state the simplest version of cell decomposition
and then formulate stronger versions as addenda:

\begin{thm}[Cell decomposition]\label{thm:cd:alg:skol.mix}
Suppose that $\acl$ equals $\dcl$ in $\Th(K)$, and that $\Th(K)$ is $1$-h-minimal.
Consider a $\emptyset$-definable set $X\subset  K^n$ for some $n$.
Then there exists a cell decomposition of $X$, namely, a partition of $X$ into finitely many $\emptyset$-definable cells $A_\ell$ of some depth $N>0$.
\end{thm}

\begin{example}
Here is an example of a cell decomposition for which one needs the condition that $\acl=\dcl$. Let $K$ be a field of equi-characteristic zero which contains cube roots of unity. Let $Z\subset K^2$ be the set of $(x,y)$ in $(K^\times)^2$ for which $x=y^3$. This is a 1-dimensional set, but it is not a cell. Denote by $X\subset K$ the set of non-zero cubes. This set $X$ is a union of fibres of the map $\rv: K\to \RV$, so it is a cell with cell center $0$. By using $\acl=\dcl$ we find three definable maps $f_1, f_2, f_3: X\to K$ such that for any $x\in X$, $\{f_1(x), f_2(x), f_3(x)\}$ is the set of cube roots of $x$. Then $Z$ is the union of three $(1,0)$-cells (of depth $1$), namely the sets
\[
Z_i = \{(x,y)\in K^2\mid (\rv(x-0), \rv(y-f_i(x))) \in \rv(X)\times \{0\}\}, \quad \text{ for } i=1,2,3.
\]
\end{example}

\begin{addendum}[Preparation of $\RV_{|M|}$-sets]\label{add:cd:alg:prep.mix}
On top of the assumptions from Theorem \ref{thm:cd:alg:skol.mix}, let also a $\emptyset$-definable set $P\subset X\times \RV_{|M|}^k$ be given for some $k$ and $M>0$. We consider $P$ as the function sending $x\in X$ to the fiber $P_x := \{\xi \in \RV_{|M|}^k \mid (x,\xi) \in P\}$.

Then the cells $A_\ell$ from Theorem \ref{thm:cd:alg:skol.mix} can be taken such that moreover $P$ (seen as function)
is constant on each twisted box of each cell $A_\ell$.
\end{addendum}

\begin{addendum}[Continuous functions]\label{add:cd:cont:f.mix}
On top of the assumptions from Theorem \ref{thm:cd:alg:skol.mix} (with Addendum~\ref{add:cd:alg:prep.mix}, if desired), suppose that finitely many $\emptyset$-definable functions $f_j\colon X \to K$ are given. Then the $A_\ell$ can be taken such that moreover, the restriction $f_j|_{A_\ell}$ of each function $f_j$ to each cell $A_\ell$ is continuous, and, that each component $c_i\colon \pi_{<i}(A_\ell) \to K$ of the center tuple of $A_\ell$ is continuous.
\end{addendum}

In the one-dimensional case, we can prepare the domain and the image of the functions in a compatible way:

\begin{addendum}[Compatible preparation of domain and image]\label{add:cd:alg:range.mix}
Under the assumptions of Addendum~\ref{add:cd:cont:f.mix}, if $n=1$, we may moreover impose that there exists an integer $M$ such that for each $\ell$ and each $j$, $f_j|_{A_\ell}$ is either constant or injective, $f_j(A_\ell)$ is a $\emptyset$-definable cell of depth $N$ and for every twisted box $R$ of $A_\ell$, there are twisted boxes $S\subset S'$ of $f_j(A_\ell)$ of depth $NM$ respectively $N$ such that
\[
S\subset f_j(R)\subset S'.
\]
\end{addendum}

In \cite{Verm:h-min}, this result is generalized in equi-characteristic zero to higher dimensions, where one can compatibly prepare the domain and image of a definable map $K^m\to K$.

We recall the usual notion of Lipschitz continuity.

\begin{defn}[Lipschitz continuity]\label{defn:Lipschitz.mix}
For a valued field $K$ and an element $\lambda$ in its value group $\Gamma_K^\times$, a function $f:X\subset K^n\to K^m$ is called Lipschitz continuous with Lipschitz constant $\lambda$ if for all $x$ and $x'$ in $X$ one has
$$
|f(x)-f(x')| \leq \lambda |x-x'|,
$$
where the norm of tuples is, as usual, the sup-norm.
We call such $f$ shortly $\lambda$-Lipschitz.

Call $f$ locally $\lambda$-Lipschitz, if for each $x\in X$ there is an open neighborhood $U$ of $x$  such that the restriction of $f$ to $U$ is $\lambda$-Lipschitz.
\end{defn}

\begin{addendum}[Lipschitz centers]
\label{add:cd:Lip:comp.mix}
Theorem~\ref{thm:cd:alg:skol.mix} (with Addenda~\ref{add:cd:alg:prep.mix} and \ref{add:cd:cont:f.mix}  if desired),
is also valid in the following variant:
Instead of imposing that $A_\ell$ itself is a cell in the sense of Definition~\ref{defn:cell.mix}, we only impose that $\sigma_{\ell}(A_\ell)$ is a cell, for some coordinate permutation $\sigma_\ell\colon K^n \to K^n$.
With this extra freedom given by the choice of $\sigma_\ell$, we can moreover ensure that there is an integer  $M>0$ such that $\sigma_{\ell}(A_\ell)$ is of cell-type $(1, \dots, 1, 0, \dots, 0)$ and that each component $c_i$ of the center tuple $(c_i)_i$ of $\sigma_{\ell}(A_\ell)$ is $|1/M|$-Lipschitz.
\end{addendum}

Closely related to that addendum, we also have the following reformulation of the piecewise continuity result of \cite{CFL} in the Hensel minimal setting.

\begin{thm}[Piecewise Lipschitz continuity]\label{thm:cd:alg:piece:Lipschitz.mix}
Suppose that $\acl$ equals $\dcl$ in $\Th(K)$ 
and that $\Th(K)$ is $1$-h-minimal.
Consider a $\emptyset$-definable set  $X\subset  K^n$ for some $n$ and a $\emptyset$-definable function $f:X\to K$. Suppose that $f$ is locally $1$-Lipschitz.
Then there exist an integer $M>0$ and a finite partition of $X$ into $\emptyset$-definable sets $A_\ell$ such that the restriction of $f$ to $A_\ell$ is $1/|M|$-Lipschitz, for each $\ell$.
\end{thm}

\begin{defn}[Supremum Jacobian property]\label{defn:sup-prep.mix}
Let integers $m>0$ and $n\ge 0$ be given.
For $X \subset K^n$ open and $f\colon X \to K$ a function, we say that $f$ has the
\emph{$|m|$-supremum Jacobian Property} 
 on $X$ if
$f$ is $C^1$ on $X$, and $f$ is either constant on $X$, or, 
for every $x$ and $y$ in $X$ with $x\not=y$ we have:
\begin{equation}\label{eq:rvn.cst}
|(\grad f)(x) - (\grad f)(y)| < |m\cdot (\grad f)(y)|
\end{equation}
and
\begin{equation}\label{eq:T3/2,mv}
|f(x) - f(y) -  (\grad f)(y)\cdot(x - y) | <  |m\cdot (\grad f)(y) |\cdot |x-y|.
\end{equation}
\end{defn}
In the definition, as usual, one considers $(\grad f)(y)$ as a matrix with a single row with entries $(\partial f/\partial x_i)(y )$, which is multiplied with the column vector $x - y$.

\begin{thm}[Sup-Jac-preparation]\label{thm:T3/2.mv.mix}
Suppose that $\Th(K)$ is $1$-h-minimal. Let integers $m>0$ and $n\ge 0$ be given. For every $\emptyset$-definable function $f\colon K^n \to K$, there exists a $\emptyset$-definable map $\chi\colon K^n\to \RV_{|N|}^k$ (for some integers $k\geq 0$ and $N\ge 1$) such that for each fiber $F$ of $\chi$ which is $n$-dimensional, $F$ is open and  $f$ restricted to $F$ has the $|m|$-supremum Jacobian property on $F$.
\end{thm}

To prove these results we first prove the cell decomposition theorem together with Addendum~\ref{add:cd:alg:prep.mix} directly. While it is possible to prove the other addenda in a similar way, we give an alternative recipe for transferring cell decompositions via coarsening.

\begin{proof}[Proof of Theorem~\ref{thm:cd:alg:skol.mix} with Addendum~\ref{add:cd:alg:prep.mix}]
We induct on $n$. For the base case $n=1$, we use Proposition~\ref{prop: 1hmixmin preparing families} to find a finite $\emptyset$-definable set $C$ and an integer $N$ such that $C$ $|N|$-prepares $P$. By our assumption that $\acl=\dcl$ every element of $C$ is $\emptyset$-definable. For each ball $B$ $|N|$-next to $C$ we can definable pick $c(B)\in C$ such that $B$ is $|N|$-next to $c(B)$. The desired cell decomposition consists of $2$ cells for every element $c\in C$, namely the $0$-cell $\{c\}$, and the $1$-cell consisting the union of all balls $B$ $|N|$-next to $C$ for which $c(B) = c$.

Now assume that $n>1$. For every $a\in K^{n-1}$, we apply the case $n=1$ to the set
\[
P_a = \{(b, \xi)\in K\times \RV_{|M|}^k\mid (a,b,\xi)\in P\}.
\]
Using compactness, we find for each $a\in K^{n-1}$ a finite set of cells of depth $|N|$ adapted to $P_a$, where we can moreover assume that the cell centers are $\emptyset$-definable functions $K^{n-1}\to K$. Now consider the set $P'\subset K^{n-1}\times \RV_{|N|} \times \RV_{|M|}^k$ where for each $a\in K^{n-1}$, $P_a'$ consists of the following:
\begin{enumerate}
\item For each $0$-cell $\{c\}\subset K$ of $P_a$, the set $\{(0, \ldots, 0)\in \RV_{|N|}\}\times P_{a,c}$ is in $P_a'$.
\item For each $1$-cell $X$ with center $c\in K$ of $P_a$ let $R\subset \RV_{|N|}$ be the set from Definition~\ref{defn:cell.mix}. Then the set
\[
\{ (\zeta, \xi)\in \RV_{|N|}\times \RV_{|M|}^k\mid \zeta \in R, \xi\in P_{a,b} \text{ for any } b\in K \text{ with } \rv_{|N|}(b-c) = \zeta\}\}
\]
is in $P_a'$.
\end{enumerate}
We find the desired cell decomposition by applying the induction hypothesis to the set $P'$.
\end{proof}

The following lemma can be used to move from a cell decomposition in $K_{\eqc}$ to one in $K$.

\begin{lem}\label{lem:celldecomp.Kecc.to.K}
Assume that $\acl$ equals $\dcl$ in $K$ and that $\Th_\cL(K)$ is 1-h-minimal. Let $\cA'$ be an $\cL_\eqc$-definable cell decomposition of $K^n$. Then there exists an $\cL$-definable cell decomposition $\cA$ of $K^n$ of depth $|N|$, for some $N$, such that every $|1|_\eqc$-twisted box of $\cA$ is contained in a twisted box of $\cA'$.
\end{lem}

By a $|1|_\eqc$-twisted box of $\cA$ we mean the following. Let $(c_1, \ldots, c_n)$ be the cell center of some cell in $\cA$ and let
\[
\{x\in K^n\mid (\rv_{|N|}(x_i-c_i(\pi_{<i}(x))))_{i=1}^n = r\}
\]
be a twisted box in this cell, for some $r\in \RV_{|N|}^n$. Then a $|1|_\eqc$-twisted box of $\cA$ is a set of the form
\[
\{x\in K^n\mid (\rv_{\eqc}(x_i-c_i(\pi_{<i}(x))))_{i=1}^n = r'\},
\]
for any $r'\in \RV_\eqc^n$ reducing to $r$ under the natural projection $\RV_\eqc\to \RV_{|N|}$.

\begin{proof}[Proof of \Cref{lem:celldecomp.Kecc.to.K}] 
We can do this one cell at a time, so let $A$ be an $\cL_\eqc$-definable cell with cell center $(c_1, \ldots, c_n)$. By \cite[Corollary 4.1.17]{CHR} (which holds in our context by \Cref{prop:resp}), we have that $c_i(\pi_{<i}(x))$ is always $\cL(\pi_{<i}(x))$-definable. Let $C_i(\pi_{<i}(x))$ be an $\cL$-definable family of finite sets, with parameter $\pi_{<i}(x)\in K^{i-1}$, such that $c_i(\pi_{<i}(x))$ is an element of $C_i(\pi_{<i}(x))$. Using $\acl=\dcl$ in $\cL$, we can find several $\cL$-definable cells of depth $|N|$, for some $N$, whose cell center is among the $C_i$. This will be the desired cell decomposition.
\end{proof}

\begin{proof}[Proof of \Cref{thm:cd:alg:piece:Lipschitz.mix,thm:T3/2.mv.mix} and the Addenda~\ref{add:cd:cont:f.mix}, ~\ref{add:cd:alg:range.mix}, ~\ref{add:cd:Lip:comp.mix}]
We explain the use of compactness to prove \Cref{thm:T3/2.mv.mix}. Let $f$ be given as in \Cref{thm:T3/2.mv.mix}. We may suppose that $K$ is saturated when we prove existence of $\chi$. Let $\cL_\eqc$ be the $\RV$-enrichment of $\cL$ obtained by adding to $\Gamma_K$ the cut given by the maximal ideal of $\cO_{K,\eqc}$. Working in $\cL_\eqc$ and applying 
\cite[Theorem~5.4.10]{CHR}, we obtain an $\cL_\eqc$-definable map $\chi': K^n\to \RV_\eqc^k$ such that on the fibres of $\chi'$, $f$ has the supremum Jacobian property with respect to $|\cdot|_\eqc$. We can assume that $\chi'$ comes from a $\cL_\eqc$-definable cell decomposition. Let $\cA$ be the $\cL$-definable cell decomposition of depth $|N|$ of $K^n$ obtained from Lemma~\ref{lem:celldecomp.Kecc.to.K} on $\chi'$. We claim that after enlarging $N$, $f$ has the $|m|$-supremum Jacobian property on the twisted boxes of $\cA$. Suppose the claim is false. By saturation, we can find $x,y\in K^n$ which are in the same $|N|$-twisted box of $\cA$ for every integer $N$ and such that either
\[
|(\grad f)(x) - (\grad f)(y)| \geq |m\cdot (\grad f)(y)|,
\]
or
\[
|f(x) - f(y) -  (\grad f)(y)\cdot(x - y) | \geq  |m\cdot (\grad f)(y) |\cdot |x-y|.
\]
However, if $x$ and $y$ are in the same $|N|$-twisted box of $\cA$ for every integer $N$, then they are also in the same $|1|_\eqc$-twisted box of $\cA$, and hence also in the same fibre of $\chi'$. But this contradicts our choice of $\chi'$, proving the claim.

The proofs of the addenda and \ref{thm:cd:alg:piece:Lipschitz.mix} can be proven in the same way using coarsening, or alternatively by direct adaptation of the proofs of \cite[Theorems 5.2.4 and 5.2.8]{CHR}.
\end{proof}


\section{Diophantine application for 1-h-minimal transcendental curves}\label{sec:counting}

In this section we give a diophantine application of Hensel minimality by bounding the number of rational points of bounded height on transcendental curves that are definable in Hensel minimal structures,
in the spirit of \cite{Bombieri-Pila}, \cite{PW}, \cite{CCL-PW} and \cite{CFL}.
In more detail, our Diophantine application for curves resembles the results for curves from Bombieri-Pila \cite{Bombieri-Pila}, but allow for an axiomatic setting, analogous to the results by Pila-Wilkie \cite{PW} (where o-minimality is used). It extends \cite{CCL-PW} and \cite{CFL} to our axiomatic setting of Hensel minimality, but only in the curve case. In \cite{CCL-PW} and \cite{CFL} the analytic (and subanalytic) situation is treated. These approaches all use Taylor approximation to finite order. In the one-dimensional case, \Cref{cor:high-ord:mix} is strong enough as Taylor approximation, but in higher dimensions new and extra work is needed which seems quite hard at the moment.  See \Cref{sec:QTr} for a discussion on the difficulties to extend the axiomatic approach to higher dimensions.

Note that the finiteness result of \cite{BCN} for the number of polynomials of bounded degree on the transcendental part of non-archimedean analytic and subanalytic sets $X$ in $K^n$ with $K=\CC\llp t \rrp$ no longer needs to hold in general $1$-h-minimal structures on $K$, and is specific to the subanalytic situation.  See \cite{CNV} for such examples with infinitely many polynomials of bounded degree on the transcendental part of a set definable in a Hensel minimal structure.
We leave the discovery of higher dimensional variants of Theorem \ref{thm:QQp-c-dim} to the future.

\subsection{}

We first define sets of rational points of bounded height, as in \cite{CCL-PW}, \cite{CFL}.

\begin{defn}
\label{defn:ratP}
For any integer $H > 0$ and any set $X\subset K^n$ with $K$ a field extension of $\QQ_p$, write $X(H)$ for 
$$
X(H) := \{x\in  X  \cap \ZZ^n\mid   0\le x_i \le H \mbox{ for each $i$} \}.
$$
\end{defn}

In the real setting, convexity together with bounds on $C^r$-norms leads to Taylor approximation results, but in our setting we have to control Taylor approximation in other ways (by lack of convexity and by total disconnectedness).
This is captured by the $T_r$ notion from \cite{CCL-PW} and \cite{CFL}.
\begin{defn}[$T_r$ maps]\label{defn:Tr}
Fix a positive integer $r$. A function $f : U\subseteq K^m \to K$ is said to be $T_r$ if for each $y\in U$ there is a polynomial $T^{<r}_{f,y}(x)$ of degree less than $r$ and coefficients in $K$ such that for each $x,y\in U$
\begin{equation}\label{eq:Tr}
\abs{f(x)-T^{<r}_{f,y}(x)}\leq \abs{x-y}^r.
\end{equation}
A map $f : U\subseteq K^m \to K^n$ is called $T_r$ if each component function is.
\end{defn}

We need some notation from~\cite[Section 5.1]{CFL}.

\begin{defn}\label{def: Notations}
Let $\alpha=(\alpha_1,...,\alpha_m)\in \NN^m$ and define $|\alpha|=\alpha_1+\cdots +\alpha_m$. We define the following sets and numbers:
\begin{itemize}
\item $\Lambda_m(k):=\{\alpha\in \NN^m;\; |\alpha|=k\}$,
\item $\Delta_m(k):=\{\alpha\in \NN^m;\; |\alpha|\leq k\}$,
\item $L_m(k):=\#\Lambda_m(k)$ and $D_m(k):=\#\Delta_m(k)$
\end{itemize}
\end{defn}

Note that $L_m(k)$ (resp.\ $D_m(k)$) is the number of monomials of degree exactly (resp.\ at most) $k$ in $m$ variables.

Fix an integer $d$ and define, for all integers $n$ and $m$ such that $m<n$ the following integers.

\begin{equation}\label{eq: def mu r V e}
\begin{aligned}
\mu(n,d) &=D_n(d) \qquad& r(m,d)&=\min\{x\in \ZZ; \; D_m(x-1)\leq \mu < D_m(x)\}\\
V(n,d) &=\sum_{k=0}^d k L_n(k) \qquad & e(n,m,d) &=\sum_{k=1}^{r-1} kL_m(k)+r(\mu-D_m(r-1)).
\end{aligned}
\end{equation}

We recall {\cite[Lemma 5.1.3]{CFL}}, amending it with a correction for the exponent (See Remark \ref{rem:amend}).

\begin{lem}[Counting algebraic hypersurfaces]
\label{lem-Tr-param-hypersurfaces}
For all positive integers $d$, $n$, $m$ with $m<n$, consider the integers $r$, $V$, $e$ and $\mu$ as defined above. Fix a finite field extension $K$ of degree $\nu$ over $\QQ_p$ and with $q_K$ many elements in the residue field of $K$. Let a subset $U\subseteq \cO_K^m$ and a $T_r$ map $\psi: U\to \cO_K^n$ be given and put $X = \psi(U)$.
Then for each $H\ge 1$, the set $X(H)$ is contained in at most
\begin{equation}\label{corr.5.1.3}
q_K^m (\mu !)^{m/e} H^{\nu mV/e}
\end{equation}
hypersurfaces of degree at most $d$.
Moreover, when $d$ goes to infinity, $mV/e$ goes to $0$.
\end{lem}

\begin{remark}\label{rem:amend}
Note that the exponent $\nu$ of $H$ in (\ref{corr.5.1.3}) is forgotten in the statement of \cite[Lemma 5.1.3]{CFL} in the characteristic zero case, but comes up in the proof given there. This lemma is used for the characteristic zero cases of Proposition 5.1.4 and Theorems B and 5.2.2 of \cite{CFL}, each of which can be corrected by letting the implicit constants depend furthermore on $\nu$. The positive characteristic cases of Lemma 5.1.3, Proposition 5.1.4 and Theorems B and 5.2.2 of \cite{CFL} don't need any change.
A sharper variant of \cite[Lemma 5.1.3]{CFL} (under some extra conditions) that would recover the original Theorems B and 5.2.2 from \cite{CFL} may be found in the future.
\end{remark}

The following lemma plays a similar role as Lemma 3.2.14 of \cite{CCL-PW} and Lemma 3.4.3 of \cite{CFL}. Similarly as in the Yomdin-Gromov and Pila-Wilkie approaches \cite{YY2} \cite{gromov} \cite{PW}, it states that composing with well-chosen powers maps yields good Taylor approximation.

\begin{lem}[Composing with power maps]\label{lem-power-maps}
Let an integer $r>0$ and a finite field extension $K$ of $\QQ_p$ be given.
Let $B$ be an open  ball of the form $d(1+N\cM_K)$ for some integer $N>0$ and some $d\in\cO_K$.  Let
$f:B\to \cO_K$ be a $T_1$ function such that  $f$ and
each of the iterated derivatives $f',f'',\ldots f^{(r)}$ of $f$ up to order $r$ has the Jacobian property, and, such that $f$ satisfies Taylor approximation up to level $r$ as in \Cref{cor:high-ord:mix}.
Then there is an integer $r'\ge r$ such that the map
$$
x\mapsto f( a+bx^{r'})
$$
 is $T_r$ on $B'$
for any choice of $a,b\in\cO_K$ and any ball $B'\subset  \cO_K$ which maps into $B$ under the map $x\mapsto a+bx^{r'}$.
\end{lem}
The proof of \Cref{lem-power-maps} is similar to the proofs of Lemma 3.4.3 of \cite{CFL} and Lemma  3.2.14 of \cite{CCL-PW} where the analyticity condition is replaced by the Jacobian property for the functions $f,f',f'',\ldots f^{(r)}$. Slightly more details of the proof can be found in \cite{CCL-PW,CFL}.
\begin{proof}[Proof of \Cref{lem-power-maps}]
The Jacobian property for $f$ and the $f^{(i)}$ and the fact that $f$ is $T_1$ imply that $|f^{(i)}(b+bNx)| \le |x|^{1-i}$ for all $i$ with $0<i\le r$ and all $x$ with $b+bNx\in B$.
Write $f_{r',B'}$ for a composition of the form $x\in B' \mapsto f( a+bx^{r'})$ for some choice of $r'$, $a,b$ and $B'$ as in the lemma.
By the chain rule for differentiation, we find that $|f_{r',B'}^{(i)}(x)| \le |r'/N^r|$ for all $i$ with $0<i\le r$ and all $x$ in $B'$. Since Taylor approximation up to level $r$ as in \Cref{cor:high-ord:mix} still holds for $f_{r',B'}$ since it is preserved under compositions (see Lemma 3.2.7 of \cite{CCL-PW}), we are done for suitable choice of $r'$.
\end{proof}

Call an infinite set $C\subset K^n$ purely transcendental if every algebraic curve in $K^n$ has finite intersection with $C$.

\begin{thm}\label{thm:QQp-c-dim}
Let $K$ be a finite field extension of $\QQ_p$.
Let $C\subset K^n$ be a purely transcendental set which is definable in a $1$-h-minimal structure on $K$
and which is of dimension $1$ (i.e.~there is a linear projection $p:K^n\to K$ such that $p(C)$ is infinite and such that $p$ has finite fibers on $C$).
Then for each $\varepsilon>0$ there is a constant $c$ such that for all $H\ge 1$ one has
$$
\# C(H) \le c    H^\varepsilon.
$$
Moreover, the constant $c$ can be chosen uniformly in a definable family of such curves in the following sense: Consider a $1$-h-minimal structure on $K$ and a definable family $\{C_y\}_{y\in Y}$ of definable sets $C_y\subset K^n$ where $y$ runs over some definable set $Y\subset K^m$.
Then for each $\varepsilon>0$ there is a constant $c$ such that for all $y\in Y$ and all $H\ge 1$ one has, if $C_y$ is of dimension $1$ and purely transcendental, then
$$
\# C_y(H) \le c    H^\varepsilon.
$$
\end{thm}

\begin{proof}
Let $K$ and $C$ be given as in the beginning of the theorem and fix $\varepsilon>0$. Write $n'$ for $\binom{n}{2}$. Fix $d$ (and corresponding value for $r$, $V$ and $e$) such that $\nu V/e< \varepsilon/n'$, with notation from \Cref{lem-Tr-param-hypersurfaces} for plane curves (namely the case of that lemma with $m=1$ and $n=2$). We may suppose that $\acl=\dcl$ by \Cref{prop:acl}, and thus,  we can speak more easily about cells. Clearly we may suppose that $C\subset \cO_K^n$.
By the Cell Decomposition Theorem \ref{thm:cd:alg:skol.mix} with Addendum \ref{add:cd:Lip:comp.mix} and possibly after composing with $x\to a+Mx$ for some $a\in \cO_K$ to kill the factor $|1/M|$ from the addendum, we may suppose that we have $N$ maps
$$
\varphi_i : P_i\subset \cO_K\to C
$$
the union of whose images equals $C$, and which are $1$-Lipschitz.
Indeed, by varying over finite choices of $a$, we can ensure that the union of images equals $C$.
Furthermore, after translation by the center we may suppose that each $P_i$ is a cell with center $0$, that each of the component functions of each partial derivative of $\varphi_i$ up to order $r$
has the Jacobian property, and, that each  component functions of $\varphi_i$ satisfies Taylor approximation up to level $r$, as in \Cref{cor:high-ord:mix}. By Lemma \ref{lem-power-maps}, we find an integer $r'\ge r$ such that the map  $x\mapsto \varphi_i( a+bx^{r'})$ is $T_r$ on $P_i'$,
for any choice of $a,b\in\cO_K$ and any cell $P_i'$ which is mapped into $P_i$ under $x\mapsto a+bx^{r'}$. Indeed, on any twisted box $B'$ in $P_i'$ this is exactly given by Lemma \ref{lem-power-maps}, and, obtaining \Cref{eq:Tr} on $P_i'$ for  $x$ and  $y$ lying in different twisted boxes of $P_i'$ is more easy.
By varying over $a$, $b$ and $P_i'$ we can ensure that the union of the images of the $\varphi_i$ still equals $C$.

Hence, we may suppose that the maps $\varphi_i$ are already $T_r$ themselves.
Write $C_i$ for $\varphi_i(P_i)$.
Apply \Cref{lem-Tr-param-hypersurfaces} to the sets  $p_j(C_i)$ in $K^2$ for each choice of coordinate projection $p_j:K^n\to K^2$. This yields $N'>0$ and for each $i,j$ and  $H\ge 1$ a collection of at most  $N' H^{\varepsilon/n'}$ algebraic curves $S_{ij\ell}$ of degree $d$ whose union contains $p_j(C_i)(H)$. Since any intersection of the form $\bigcap_j p_j^{-1}(S_{ij\ell})$ is algebraic and of dimension at most $1$ (see for example the argument on transversal cylinders  on page 45 of \cite{CCL-PW}), 
we have for each $H\ge 1$ that the set $C_i(H)$ is contained in no more than $N'{}^{n'} H^{\varepsilon}$ many algebraic curves (all of bounded complexity in terms of $d$ and $n$). Each such algebraic curve has finite intersection with $C_i$ by the pure transcendence of $C_i$. Since the intersections of  $C_i$ with the mentioned algebraic curves appear in a definable family of finite definable sets, and
since one has uniform upper bounds on the size of finite sets in definable families in $1$-h-minimal structures by (3) of \Cref{lem: basics 1hmixmin}, the theorem follows.
\end{proof}

\begin{remark}
A version of \Cref{thm:QQp-c-dim} which works uniformly in all local fields $K$ of large residue field characteristic (namely $\QQ_p$ and $\FF_p\llp t\rrp$ and their finite field extensions, for large $p$) can also be formulated and proved along the same lines. This extends the one-dimensional case (with $m=1$, and with implied constant as specified in Remark \ref{rem:amend}) of Theorem 5.2.2 of \cite{CFL} to the $1$-h-minimal situation (instead of subanalytic).
\end{remark}

\section{Some questions}

We end the paper with some questions for future research.


\subsection{}

Suppose that $\cL'$ is an \(\bigcup_n\RV_n\)-expansion of $\cL$. Then $\ell$-\hmix-minimality of $\Th_{\cL}(K)$ implies $\ell$-\hmix-minimality of $\Th_{\cL'}(K)$
for $\ell = 0, 1, \omega$ (see \Cref{prop: RVenrichment mixed} and \Cref{thm:equivalences:1hmin}).

\begin{question}\label{que:conv}
Does the converse also hold, i.e., does $\ell$-h-minimality of $\Th_{\cL'}(K)$ imply $\ell$-h-minimality of $\Th_{\cL}(K)$?
\end{question}

\begin{remark}
Hensel minimality is not preserved by passing to reducts in general. Indeed, suppose that $\Th_{\cL}(K)$ is $\omega$-h-minimal and that $K$ is $\aleph_0$-saturated.
Fix a ball $B = B_{<\lambda}(a) \subset K$ which is strictly contained in a ball disjoint from $\acl_K(\emptyset)$
(so that $B$ cannot be prepared by a finite, $\emptyset$-definable $C \subset K$). Then $\Th_{\cL(a,\lambda)}(K)$ is $\omega$-h-minimal but
the reduct $\Th_{\cL \cup \{B\}}(K)$ is not even $0$-h-minimal (where by ``$B$'', we mean a predicate for that ball). It would be interesting to find a tameness notion which is preserved on reducts.
\end{remark}

\subsection{}

Suppose that $\Th(K)$ is $\omega$-\heqc-minimal and that we have a definable coarsening $|\cdot|_c$ of the valuation;
write $K_c$ for $K$ considered as a valued field with the coarsened valuation $|\cdot|_c$, and write $k_c$ for the residue field of $K_c$,
and, put the full induced structure on $k_c$.
By 
\cite[Corollary~4.2.4]{CHR}, resp.~by \Cref{thm:equivalences:1hmin}, one has
\begin{enumerate}
\item if $\Th(K)$ is $\omega$-\heqc-minimal, then so is $\Th(K_c)$, and, resp., if $\Th(K)$ is $1$-h-minimal, then so is $\Th(K_c)$.
\end{enumerate}

\begin{question}
\begin{enumerate}\stepcounter{enumi}
 \item Does $\omega$-\heqc-minimality of $\Th(K)$ imply $\omega$-\heqc-minimality of $\Th(k_c)$?
 \item Do $\omega$-\heqc-minimality of $\Th(K_c)$ and $\omega$-\heqc-minimality of $\Th(k_c)$ together imply $\omega$-\heqc-minimality of $K$?
\end{enumerate}
And one may ask (2) and (3) also for $1$-h-minimality instead of $\omega$-\heqc-minimality. One may also ask whether  (1) holds under other variants of Hensel minimality.
\end{question}

\subsection{}
Any $C^1$-function $U \subset \RR^n \to \RR$ also has a strict derivative (see 
\cite[Definition~3.1.2]{CHR}.
This is not the case in valued fields:

\begin{example}
Define $f\colon K^2 \to K$ by $f(x,y)= x^2$ if $|x|^4 \leq  |y|$ and $f(x,y) = x^3$ otherwise. This function is $C^1$ everywhere, but at $0$, the strict derivative does not exist, since $\frac{f(x,x^4) - f(x,0)}{x^4} = x^{-2} - x^{-1}$, which diverges for $x \to 0$.
\end{example}

In view of this example, and by our knowledge that strict $C^1$ is the better notion for rank one valued fields for several reasons (see e.g.~\cite{BGlockN}, and where strict $C^1$ means that the strict derivative exists everywhere and is continuous), one may try to build a good working notion of definable strict $C^1$ submanifolds of $K^n$, assuming a suitable form of Hensel minimality.   The following is a first question in this direction.

\begin{question}
In what form could an Implicit Function Theorem hold for definable functions, say, assuming $1$-h-minimality? This should be considered knowing that a $1$-h-minimal theory does not automatically imply Hensel's lemma in mixed characteristic, and, knowing that, in some special situations (like in \cite{BGlockN}) it is known that strict $C^1$ functions work well for the Implicit Function Theorem.
\end{question}

%
%


\subsection{}\label{sec:QTr}

\Cref{thm:T3/2.mv.mix} can be seen as a order one Taylor approximations result, for functions in several variables. This viewpoint suggests
that we might have the following variant of 
\cite[Theorem~5.6.1]{CHR} (which is a result on higher order Taylor approximations of functions in several variables):

\begin{question}\label{qu:t-high-high}
Given a definable function $f\colon K^n \to K$ in a $1$-h-minimal structure and an integer $r \ge 1$, does there exist a definable map $\chi\colon K^n \to \RV_{|N|}^k$ such that
\cite[Equation~(5.6.1)]{CHR} (or a similar kind of Taylor approximation) holds on each $n$-dimensional fiber of $\chi$?
\end{question}

Such a result would be strictly stronger than 
\cite[Theorem~5.6.1]{CHR}, which yields Taylor approximations
only on boxes disjoint from a lower-dimensional definable set $C$. Indeed, given $\chi$, one can easily find a $C$
such that every box disjoint from $C$ is contained in a fiber of $\chi$
(namely by $1$-preparing $\chi$ fiberwise using \Cref{cor: 1hmixmin preparing families}). On the other hand,
the family of maximal boxes disjoint from $C$ cannot, in general, by parametrized by a tuple from $\RV$. An answer to \Cref{qu:t-high-high} is important to generalize the diophantine application of \Cref{sec:counting} to arbitrary dimension instead of just curves.

\bibliographystyle{amsplain}
\bibliography{anbib}
\end{document}